\documentclass[12pt]{article}
\usepackage{epsf,epsfig,color}
\usepackage{amssymb,latexsym}
\usepackage{graphicx}%

\newtheorem{lem}{Lemma}[section]
\newtheorem{theo}{Theorem}[section]
\newtheorem{pro}{Proposition}[section]
\newtheorem{cor}{Corollary}[section]

\newtheorem{prob}{Problem}[section]

\renewcommand{\theenumi}{\rm (\roman{enumi})}

\newcommand{\proof}
{{\noindent {\em Proof}.\quad}\setcounter{countclaim}{0}
\setcounter{countcase}{0}}
\newcommand{\proofend}{{\hfill$\Box$}}

\newcounter{countcase}

\newcounter{countclaim}
\def\inclaim{\addtocounter{countclaim}{1}
{\noindent {\bf Claim \thecountclaim}: }}

\newcounter{countfig}

\newcommand \spann[1]
{\langle #1\rangle}

\newcommand{\beeq}{\begin{equation}}
\newcommand{\eneq}{\end{equation}}

\newcommand{\beeqn}{\begin{eqnarray*}}
\newcommand{\eneqn}{\end{eqnarray*}}

\def \sets{{\cal S}}

\def\setn{{\cal N}}

\def \setst{{\cal ST}}

\def \R {{\cal R}}

\def \N {{\mathbb N}}

\def \iff {if and only if }

\newcommand {\relabel}[1] {\label{#1} \red{[*: #1]}}
\newcommand {\rebibitem}[1] {\bibitem{#1} \red{[*: #1]}}
\def\relabel {\label} 
\def\rebibitem {\bibitem}  

\begin{document}

\newcommand{\resection}[1]
{\section{#1}\setcounter{equation}{0}}

\renewcommand{\theequation}{\thesection.\arabic{equation}}

\renewcommand{\theenumi}{\rm (\roman{enumi})}
\renewcommand{\labelenumi}{\rm(\roman{enumi})}

\baselineskip 0.6 cm

\title {Formulas counting spanning trees in line graphs and their extensions
}

\author
{Fengming Dong\thanks{Corresponding author.
Email: fengming.dong@nie.edu.sg}\\
\small National Institute of Education\\
\small Nanyang Technological University, Singapore
}

\date{}

\maketitle

\begin{abstract}
For any connected multigraph $G=(V,E)$ and any $M\subseteq E$,
if $M$ induces an acyclic subgraph of $G$
and removing all edges in $M$ yields a subgraph of $G$ 
whose components are complete graphs,
a formula for $\tau_G(M)$ is obtained,
where $\tau_G(M)$  is the number of spanning trees in $G$ 
which contain all edges in $M$. 
Applying this result, we can easily obtain 
a formula for the number of spanning trees in 
the line graph or the middle graph of an arbitrary 
graph.
Applying this result, we also show that 
for any connected graph $G$ with a clique $U$ 
which is a cut-set of $G$, 
the number of spanning trees in $G$ 
has a factorization 
which is analogous to a property 
of the chromatic polynomial of $G$.
\end{abstract}

\noindent {\bf MSC}: 05A15, 05C05, 05C30, 05C76

\noindent {\bf Keywords}: Graph, clique,
Spanning tree, Cayley's formula

\resection{Introduction}

The graphs considered in this article 
are multigraphs without loops.
For any graph $G$,
let $V(G)$ and $E(G)$ 
be the vertex set and the edge set  
of $G$ respectively.
For any non-empty $V'\subseteq V(G)$, 
let $G[V']$ denote the subgraph of $G$ induced by $V'$,
and when $V'\ne V(G)$, 
let $G-V'$ be the subgraph $G[V(G)-V']$ 
(i.e., the subgraph of $G$ obtained by deleting all 
vertices in $V'$).
Let $N_G(V')=\bigcup_{v\in V'}N_G(v)$,
where $N_G(v)$ is the set of neighbours of $v$ in $G$,
and $N_G[V']=V'\cup N_G(V')$. 
For any $E'\subseteq E(G)$, 
let $G\spann {E'}$ be the spanning subgraph of $G$ 
with edge set $E'$,
let $G[E']$ be the subgraph of $G$ induced by $E'$ 
when $E'\ne \emptyset$ (i.e., the graph obtained from 
$G\spann {E'}$ by removing all isolated vertices),
let $G/E'$ be the graph obtained from $G$ by contracting 
all edges in $E'$
and $G-E'$ be the subgraph $G\spann{E(G)-E'}$
(i.e., the graph obtained from $G$ by removing all edges in $E'$).

For any graph $G$, let $\setst_G$ be the set of spanning trees 
of $G$ and let $\tau_G=|\setst_G|$.
Clearly, $\tau_G=0$  \iff $G$ is disconnected. 
It is well-known that 
$\tau(K_n)=n^{n-2}$, due to Cayley~\cite{Aig},
where $K_n$ is the complete graph of order $n$.
This beautiful formula was extended
by Moon \cite{moon0,moon1,moon2}
(also see Lov\'{a}sz  \cite[Problem 4 in page 34]{lov})
for counting the number of 
spanning trees $T\in \setst_{K_n}$ 
which contain all edges of a given forest in $K_n$.

For any $M\subseteq E(G)$, 
let $\setst_G(M)$ be the set of those members
$T\in \setst_G$ with $M\subseteq E(T)$
and let $\tau_G(M)=|\setst_G(M)|$.
Thus $\setst_G(M)\subseteq \setst_G$,
where $\setst_G(M)=\setst_G$ holds whenever 
$M$ consists of bridges of $G$.
Clearly, 
$\tau_G(M)=0$ \iff 
either $G$ is disconnected or $G\spann {M}$ contains cycles. 

\begin{theo}[Lov\'{a}sz~\cite{lov} and Moon~\cite{moon0,moon1,moon2}]\relabel{th1-1}
For any $M\subseteq E(K_n)$, 
if $K_n\spann {M}$ is a forest with $c$ components 
whose orders are $n_1,n_2,\cdots,n_c$,  then 
\begin{equation}\relabel{eq1-1}
\tau_{K_n}(M)=n^{c-2}\prod_{i=1}^c n_i.
\end{equation}
\end{theo}

It is natural to consider a suitable extension of 
Theorem~\ref{th1-1}.
In this article, 
we assume that 
$G=(V,E)$ is a connected graph, where  
$V$ can be partitioned into subsets $V_0, V_1,\cdots,V_k$
and $V_i$ is a clique of $G$ 
(i.e., $G[V_i]$ is a complete graph)
for all $i=1,2,\cdots,k$.
Thus $G[V_i]$ has no parallel edges for all $1\le i\le n$,
although $G$ may have parallel edges. 
Note that $V_0$ may be an empty set 
and $G[V_0]$ may be not complete
and may have parallel edges also. 

For any $U_1, U_2\subseteq V$, 
let $E_G(U_1,U_2)$ denote the set 
of those edges in $G$ with one end in $U_1$ and 
another end in $U_2$,
and let $E_G(U_1)=E_G(U_1, V-U_1)$.
In the case that $V_0=\emptyset$ and $M_0=\bigcup_{1\le i<j\le k}E_G(V_i,V_j)$ is a matching of $G$, 
an formula for $\tau_G(M_0)$ was obtained in 
\cite[Theorem 3.1]{dong}.
Let $G^*$ be the graph obtained from $G$ 
by identifying all vertices in each $V_i$ as one vertex $v_i$
for all $i=1,2,\cdots,k$
and removing all loops.
Thus $G^*=G/E_0$, where $E_0=\cup_{1\le i\le k}E(G[V_i])$.

\begin{theo}[\cite{dong}]\relabel{th1-2} 
If $V_0=\emptyset$ and $M_0$ is a matching of $G$,
then 
\begin{equation}\relabel{eq1-2}
\tau_G(M_0)
=\prod_{i=1}^k |V_i|^{|V_i|-2}
\sum_{T\in \setst_{G^*}} \prod_{e\in M_0-E(T)}
\left ( 
|V_{a(e)}|^{-1}+|V_{b(e)}|^{-1}
\right ),
\end{equation} 
where $1\le a(e)<b(e)\le k$ such that 
$v_{a(e)}$ and $v_{b(e)}$ are the two ends of $e$ in $G^*$
for each $e\in E(T)$.
\end{theo}

If $V_0=\emptyset$ and $M_0$ is a perfect matching of $G$, 
then $G/M_0$ is actually the line graph $L(G^*)$ of $G^*$. 
Since $\tau_G(M_0)=\tau_{G/M_0}$
holds 
(see Lemma~\ref{le2-01}~\ref{le2-1-n5}), 
applying Theorem~\ref{th1-2} yields 
a relation between $\tau_{L(H)}$
and $\tau_{H}$
for an arbitrary connected graph $H$.

\begin{cor}[\cite{dong}]\relabel{co1-1}
Let $H$ be a connected and loopless graph with vertices 
$v_1,v_2,\cdots, v_k$. 
Then 
\begin{eqnarray}\relabel{eq1-3}
\tau_{L(H)}
&=&\prod_{i=1}^k |d(v_i)|^{|d(v_i)|-2}
\sum_{T\in \setst_{H}} 
\prod_{e\in E(H)-E(T)}
\left ( 
|d(v_{a(e)})|^{-1}+|d(v_{b(e)})|^{-1}
\right ).
\end{eqnarray} 
where $d(v_i)$ is the degree of $v_i$ in $H$
and $v_{a(e)}$ and $v_{b(e)}$ are the two ends of $e$ in $H$.
\end{cor}

For any connected graph $H$, 
the middle graph $M(H)$ of $H$ 
is the one obtained from $H$ 
by subdividing each edge in $H$ exactly once
and adding a new edge joining each pair of 
new vertices $u_1,u_2$ which subdivide 
a pair of adjacent edges in $H$ (see \cite{che}).

Observe that if $V_0=\emptyset$, 
$M_0$ is a matching of $G$ and exactly one vertex in each 
$V_i$ is not incident with $M_0$, where $1\le i\le k$, 
then $G/M_0$ is actually the middle graph of $G^*$,
and thus, by the equality $\tau_G(M_0)=\tau_{G/M_0}$,
a formula for $\tau_{M(H)}$ 
follows directly from Theorem~\ref{th1-2}.

\noindent {\bf Remark}: 
The study of a relation between $\tau_{L(H)}$ 
and $\tau_H$ for a connected graph $H$ 
was started in 1966 when 
Vahovskii \cite{vah} first established such a relation
for a $r$-regular graph $H$: 
\begin{equation}\relabel{eq1-4}
\tau_{L(H)}=2^{m-n+1}r^{m-n-1}\tau_{H},
\end{equation}
where $n=|V(H)|$ and $m=|E(H)|$.
When $H$ is a graph in which each vertex is of degree 
$1$ or $r$, where $r$ is a constant, 
a similar relation between $\tau_{L(H)}$ and $\tau_H$
was found by Yan~\cite{yan1} in 2013.
When $H$ is an $(a, b)$-semiregular bipartite graph,
such a relation 
was found by Cvetkovi\'{c} (see \cite[see Theorem 3.9]{mac}, 
\cite[\S 5.2]{moh}, or \cite {sat}).
Corollary~\ref{co1-1} was the first result giving  
a relation between $\tau_{L(H)}$ and $\tau_H$ 
for an arbitrary connected graph $H$,
which implies all these known results.

In this article, we will further extend Theorem~\ref{th1-2}.
Recall that $V_0,V_1,\cdots,V_k$ 
is a partition of $V$, where 
$V_i$ is a clique of $G$ 
for all $i=1,2,\cdots,k$.
Let $M=M_0\cup \bigcup_{1\le i\le k} E_G(V_0,V_i)$,
where $M_0=\bigcup_{1\le i<j\le k}E_G(V_i,V_j)$ is not restricted to a matching of $G$. 
We will study the set $\setst_G(W)$ for any 
$W$ with $M\subseteq W\subseteq E(G)$,
where $G\spann {W}$ is a forest.

In Section~\ref{sec2}, 
we will transform $W, M$ and $G$
so that 
the study of $\tau_G(W)$ 
can be restricted to the special case 
that $M_0=\emptyset$, each component of $G[M]$ 
is a star with a center in $V_0$
and $W=M\cup N$ for some $N\subseteq E(G[V_0])$,
as stated in (i), (ii) and (iii) in Page~\pageref{special-case}.
With these conditions, 
the structure of $G$ is as shown in Figure~\ref{f4-1}
of Page~\pageref{special-case}.
Thus, in Sections~\ref{sec3} and~\ref{sec4},
we will study $\tau_G(W)$ under these assumptions.

For any $U\subseteq V(G)$, let $G\bullet U$ denote the 
graph $G/E(G[U])$, i.e., the graph obtained from $G$ by 
contracting all edges in $G[U]$. 
In Section~\ref{sec3}, 
we find a relation between $\tau_G(M\cup N)$ 
and $\tau_{G\bullet U}(N)$, 
where $U$ is a clique of $G$ and $N\subseteq E(G[V-U])$.
In Section~\ref{sec4},
we will apply the result in Section~\ref{sec3}
to obtain a relation between $\tau_G(M\cup N)$ 
and $\tau_{G\bullet U}(N)$,
where $U$ is the union of 
$k$ disjoint cliques $V_1, V_2, \cdots,V_k$
and $N\subseteq E(G[V-U])$.

In Section~\ref{sec5}, 
as an application of the results in Section~\ref{sec4},
we find a formula for $\tau_G$
when $E(G)$ can be partitioned into subsets 
$E_1,E_2,\cdots,E_k$ such that 
each $G[E_i]$ is a clique in $G$.
The middle graph and the line graph of any connected graph $H$  
are examples of such graphs.
Applying the results in Section~\ref{sec4},
one can easily deduce formulas for 
$\tau_{M(H)}$ and $\tau_{L(H)}$
for any given connected graph $H$.

It is well known that if $U$ is a clique of $G$ 
and $S_1,S_2$ is a partition of $V-U$ such that 
$E_G(S_1, S_2)=\emptyset$, then 
the following equality for the chromatic 
polynomial $\chi(G,\lambda)$ of $G$ holds 
(see~\cite{dong0,rea1, zyk}):
\begin{equation}\relabel{eq1-5}
\chi(G,\lambda)=\frac{
\chi(G[U\cup S_1], \lambda)\cdot 
\chi(G[U\cup S_2], \lambda)}
{\chi(K_{|U|},\lambda)}.
\end{equation}
Section~\ref{sec6} shows that 
$\tau_G$ has a similar result 
as (\ref{eq1-5}) when 
$N_G[S_1]\cap N_G[S_2]=\emptyset$ holds.

\resection{
It suffices to study $\tau_G(W)$ 
for a special case
\relabel{sec2}
}

Let $G=(V,E)$ be a connected graph 
whose vertex set has a partition
$V_0,V_1,\cdots,V_k$,
where each $V_i$ is a clique for all $i=1,2,\cdots,k$.
In this section, we will show that 
for any set $W$ with 
$\bigcup_{0\le i\le k}E_G(V_i,V_j)\subseteq W\subseteq E$,
the study of $\tau_G(W)$ can be 
transformed to the special case that 
each component of $G[W]$ is a star
and $E_G(V_i,V_j)=\emptyset$ for all $i,j$ 
with $1\le i<j\le k$.

\subsection{
$\tau_G(W)=\tau_{G\star W}(W')$ holds
for $W'=E(G\star W)-E(G)$
\relabel{sec2-1}
}

For any $E'\subseteq E$, 
one can easily prove 
the following basic properties on $\setst_G(E')$.

\begin{lem}\relabel{le2-01}
Let $E'\subseteq E$ and $e\in E$.
Then 
\begin{enumerate}
\item\relabel{le2-1-n1}
if $G\spann{E'}$ contains cycles, then 
$\setst_G(E')=\emptyset$;
\item\relabel{le2-1-n2} 
if $e\notin E'$ and $e$ is a loop, 
$\setst_G(E')=\setst_{G-e}(E')$;
\item\relabel{le2-1-n3}
if $e\notin E'$ and $e$ is parallel to an edge in $E'$, then $\setst_G(E')=\setst_{G-e}(E')$;
\item\relabel{le2-1-n4} if $e\in E'$ and $e$ is not a loop, 
then $\tau_G(E')=\tau_{G/e}(E'-\{e\})$;
\item\relabel{le2-1-n5} 
if $G\spann{E'}$ is a forest,
$\tau_G(E')=\tau_{G/E'}$;
\item\relabel{le2-1-no6} if $e\notin E'$
and $G[E'\cup \{e\}]$ has a cycle containing $e$, 
then $\setst_G(E')=\setst_{G-e}(E')$.
\end{enumerate}
\end{lem}

\proof \ref{le2-1-n1} and \ref{le2-1-n3}
follow directly from the definition of $\setst_G(E')$.
Both \ref{le2-1-n2} and 
\ref{le2-1-no6} follow from the fact 
that $e$ is not contained in any tree $T\in \setst_G(E')$.
\ref{le2-1-n4} follows from the fact that 
$T\in \setst_G(E')$ \iff $T/e\in \setst_{G/e}(E'-\{e\})$.
\ref{le2-1-n5} follows from \ref{le2-1-n4} directly.
\proofend

For any $W\subseteq E$,
let $G\star W$ 
denote the graph obtained from $G$ \relabel{starw}
by adding a new vertex $w_i$ and new edges 
joining $w_i$ to all vertices in $W_i$
for all $i=1,2,\cdots,r$,
where $W_1,\cdots,W_r$ are the components of $G[W]$.
An example of $G\star W$ is shown in Figure~\ref{f1},
where $G[W]$ has three components.
Thus, $V(G\star W)=V(G)\cup \{w_1,w_2,\cdots,w_r\}$ 
and $E(G\star W)=E(G)\cup \bigcup_{1\le i\le r}E_{G\star W}(w_i)$.
Also note that $\{w_1,w_2,\cdots,w_r\}$ is an independent set 
in $G\star W$.

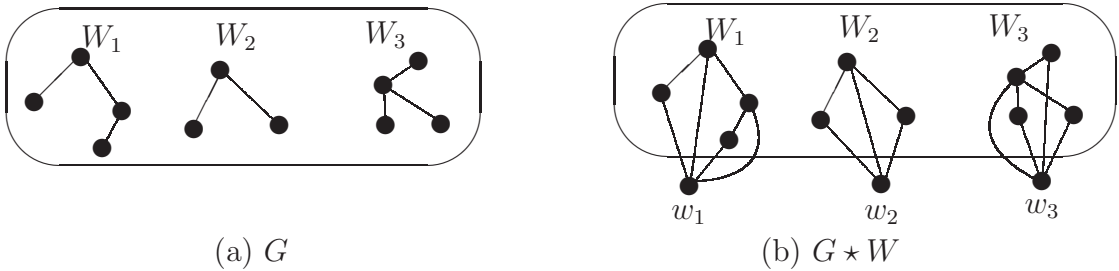
\begin{figure}[h!]
 \centering
\unitlength .9mm 
\linethickness{0.4pt}
\ifx\plotpoint\undefined\newsavebox{\plotpoint}\fi 
\begin{picture}(165.109,32.35)(0,0)
\put(5.489,18.175){\line(1,1){6.832}}
\put(98.12,19.413){\line(1,1){6.832}}
\multiput(12.322,25.007)(.03744375,-.0532125){160}{\line(0,-1){.0532125}}
\multiput(104.953,26.245)(.03744375,-.0532125){160}{\line(0,-1){.0532125}}
\multiput(18.313,16.493)(-.03725316,-.06918987){79}{\line(0,-1){.06918987}}
\multiput(110.944,17.731)(-.03725316,-.06918987){79}{\line(0,-1){.06918987}}
\multiput(33.134,22.379)(.039413462,-.037394231){208}{\line(1,0){.039413462}}
\multiput(125.765,23.617)(.039413462,-.037394231){208}{\line(1,0){.039413462}}
\multiput(56.823,20.31)(.05488889,.03737778){90}{\line(1,0){.05488889}}
\multiput(150.337,21.548)(.05488889,.03737778){90}{\line(1,0){.05488889}}
\multiput(57.033,20.205)(.059034247,-.037438356){146}{\line(1,0){.059034247}}
\multiput(150.547,21.443)(.059034247,-.037438356){146}{\line(1,0){.059034247}}
\put(15.37,11.027){\circle*{2.828}}
\put(108.001,12.265){\circle*{2.828}}
\put(102.129,5.501){\circle*{2.828}}
\put(130.447,5.724){\circle*{2.828}}
\put(154.18,6.142){\circle*{2.828}}
\put(18.313,16.493){\circle*{2.828}}
\put(110.944,17.731){\circle*{2.828}}
\put(12.216,24.481){\circle*{2.828}}
\put(104.847,25.719){\circle*{2.828}}
\put(5.384,17.859){\circle*{2.828}}
\put(98.015,19.097){\circle*{2.828}}
\put(32.818,22.589){\circle*{2.828}}
\put(125.449,23.827){\circle*{2.828}}
\put(28.928,13.865){\circle*{2.828}}
\put(121.559,15.103){\circle*{2.828}}
\put(41.543,14.496){\circle*{2.828}}
\put(134.174,15.734){\circle*{2.828}}
\put(56.928,20.31){\circle*{2.828}}
\put(150.442,21.548){\circle*{2.828}}
\put(62.078,23.884){\circle*{2.828}}
\put(155.592,25.122){\circle*{2.828}}
\put(65.442,14.634){\circle*{2.828}}
\put(158.956,15.872){\circle*{2.828}}
\put(57.243,14.529){\circle*{2.828}}
\put(150.757,15.767){\circle*{2.828}}
\multiput(130.81,5.538)(.03709412,.11747059){85}{\line(0,1){.11747059}}
\multiput(125.765,23.722)(.03737037,-.132362963){135}{\line(0,-1){.132362963}}
\multiput(97.91,18.992)(.037478261,-.118817391){115}{\line(0,-1){.118817391}}
\multiput(102.22,5.328)(.03734211,.27106579){76}{\line(0,1){.27106579}}
\multiput(107.791,12.055)(-.037438356,-.045356164){146}{\line(0,-1){.045356164}}
\qbezier(110.944,18.046)(116.094,6.589)(102.325,6.063)
\multiput(154.121,5.886)(.037424242,.080431818){132}{\line(0,1){.080431818}}
\multiput(155.277,25.227)(-.0362414,-.6741724){29}{\line(0,-1){.6741724}}
\qbezier(150.652,21.548)(140.929,12.614)(153.7,6.202)
\put(102.067,1.228){\makebox(0,0)[cc]{$w_1$}}
\put(130.552,1.123){\makebox(0,0)[cc]{$w_2$}}
\put(154.283,2.102){\makebox(0,0)[cc]{$w_3$}}
\multiput(130.687,5.504)(-.0374634146,.0398821138){246}{\line(0,1){.0398821138}}
\put(28.84,14.077){\line(1,2){4.237}}
\put(121.471,15.315){\line(1,2){4.237}}
\multiput(150.778,15.659)(.037170213,-.105957447){94}{\line(0,-1){.105957447}}
\multiput(57.264,14.346)(-.037125,.734){8}{\line(0,1){.734}}
\multiput(150.778,15.584)(-.037125,.734){8}{\line(0,1){.734}}
\put(36.239,20.064){\oval(70.003,23.158)[]}
\put(128.074,21.036){\oval(74.069,22.627)[]}
\put(15.38,27.223){\makebox(0,0)[cc]{$W_1$}}
\put(107.48,28.284){\makebox(0,0)[cc]{$W_1$}}
\put(35.178,27.577){\makebox(0,0)[cc]{$W_2$}}
\put(127.279,28.638){\makebox(0,0)[cc]{$W_2$}}
\put(57.275,27.754){\makebox(0,0)[cc]{$W_3$}}
\put(149.376,28.814){\makebox(0,0)[cc]{$W_3$}}
\end{picture}


(a) $G$ \hspace{6 cm} (b) $G\star W$

\caption{Graph $G\star W$, where $G[W]$ has 3 components}
\relabel{f1}
\end{figure}

\begin{lem}\relabel{le2-11}
Let $W\subseteq E$ such that $G[W]$ is a forest.
For any $W_0\subseteq W$, 
\begin{equation}\relabel{eq2-1}
\tau_G(W)=\tau_{G\star W}(W')
=\tau_{G\star W-W_0}(W'),
\end{equation}
where $W'=E(G\star W)-E(G)$.
\end{lem}

\proof Observe that 
$\tau_{G\star W}(W')
=\tau_{G\star W-W_0}(W')$ follows from 
Lemma~\ref{le2-01}~\ref{le2-1-no6} directly,
while $\tau_G(W)=\tau_{G\star W}(W')$ follows from 
Lemma~\ref{le2-01}~\ref{le2-1-n5} and the fact that 
the two graphs obtained respectively 
from $G/W$ and $(G\star W)/W'$
by removing their loops are isomorphic.  

Thus the result follows. 
\proofend

\subsection{Transformed to a special case}
\relabel{sec2-2}

Recall that  $V_0,V_1,\cdots,V_k$ is a partition of $V$ 
such that each $V_i$ is a clique of $G$ for all 
$i=1,2,\cdots,k$.
Let $M=\bigcup_{0\le i<j\le k}E_G(V_i,V_j)$
and $W$ be a subset of $E$ with $M\subseteq W$ 
such that $G[W]$ is a forest.
By Lemma~\ref{le2-11}, we get the following conclusion.


\begin{lem}\relabel{le2-3}
Let $G'$ denote the graph $G\star W-M$, $V'_0=V_0\cup (V(G\star W)-V(G))$ and $W'=E(G\star W)-E(G)$. 
The following properties hold:
\begin{enumerate}
\item\relabel{le2-3-n1} $\tau_G(W)=\tau_{G'}(W')$; 
\item\relabel{le2-3-n2} 
$V'_0,V_1,\cdots,V_k$ is a partition of $V(G')$,
where $V_i$ is a clique of $G'$
for all $i=1,2,\cdots,k$;
\item\relabel{le2-3-n3} $E_{G'}(V_i,V_j)=\emptyset$
for all $i,j$ with $1\le i<j\le k$;
\item\relabel{le2-3-n4} 
each component of $G'[W']$ is a star with a center in 
$V(G')-V(G)\subseteq V'_0$.
\end{enumerate}
\end{lem}

Lemma~\ref{le2-3}~\ref{le2-3-n1} follows from 
Lemma~\ref{le2-11} while 
Lemma~\ref{le2-3}~\ref{le2-3-n2}-\ref{le2-3-n4} follow directly from the definitions of $G'$ 
and $W'$.
By Lemma~\ref{le2-3},
the study of $\tau_G(W)$ can be 
restricted to the special case 
that $V(G)$ has a partition $V_0, V_1, \cdots, V_k$ 
satisfying the following conditions: 
\begin{enumerate}
\item $V_i$ is a clique for all $i=1,2,\cdots,k$
and $E_{G}(V_i,V_j)=\emptyset$ holds
for each pair of $i,j$ with $1\le i<j\le k$;
\item for $M=\bigcup_{0\le i\le k}E_G(V_0,V_i)$,
each component of $G[M]$ is a star with a center 
in $V_0$ (i.e., $d_G(u)\le |V_i|$ holds 
for each $u\in V(G)-V_0$);
\item $W=M\cup N$ for some $N\subseteq E(G[V_0])$. 
\relabel{special-case}
\end{enumerate}
When the above three conditions holds, 
$G$ has its structure as shown in Figure~\ref{f4-1}.

\begin{figure}[h!]
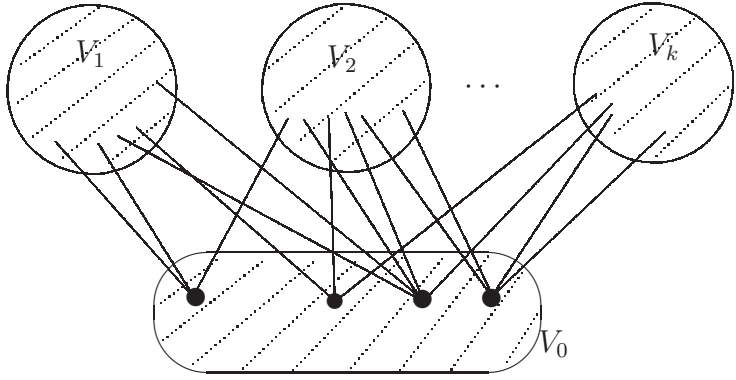

 \centering
\unitlength .8mm 
\linethickness{0.4pt}
\ifx\plotpoint\undefined\newsavebox{\plotpoint}\fi 



\caption{$E_G(V_i,V_j)=\emptyset$ for all 
$1\le i<j\le k$ and
each component of $G[M]$ is a star with a center 
in $V_0$ 
}
\relabel{f4-1}
\end{figure}

\resection{
Contracting a clique $U$ 
\relabel{sec3}}

Let $U$ be a clique of a connected graph $G=(V,E)$.
In Subsection~\ref{sec3-1}, we will deduce a formula 
for $\tau_G(W)$ in the case that
$G[W]$ is a forest,
where $W=E-E(G[U])$.
Let $G\bullet U$ denote the graph $G/G[U]$.
In Subsection~\ref{sec3-2},
we will give a relation between 
$\tau_G(M\cup N)$ and 
$\tau_{G\bullet U}(N)$,
where $M=E_G(U)$ and $N\subseteq E(G-U)$,
under the condition that 
each component of $G[M]$ is a star with a center in $V-U$.

\subsection{
When $G-E(G[U])$ is a forest
\relabel{sec3-1}}

Note that $G-E(G[U])$ is a forest \iff 
$G[W]$ is a forest, where $W=E-E(G[U])$. 
When $G[W]$ is a forest, 
applying Theorem~\ref{th1-1},
we get a formula for 
$\tau_G(W)$ below.

\begin{pro}
\relabel{pro3-1} 
Let $U$ be a clique of $G$ with $U\ne V$
and $W=E-E(G[U])$.
If $F=G[W]$ is a forest with components
$F_1,F_2,\cdots, F_t$, 
then  
\begin{equation}\relabel{eq3-20}
\tau_G(W)=
|U|^{|U|-2+t-n_1-n_2-\cdots-n_t}
\prod_{i=1}^t n_i,
\end{equation}
where $t$ is the number of components of $F$ and 
$n_i=|V(F_i)\cap U|$ for $i=1,2,\cdots,t$.
\end{pro}

\proof
Note that for any $i=1,2,\cdots,t$, 
$|E(F_i)|\ge n_i$, and 
$|E(F_i)|= n_i$ \iff $F_i$ is a star 
with a center in $V-U$ 
and $E(F_i)\subseteq E_G(U)$.

We shall prove this result by the following claims.

\inclaim  
(\ref{eq3-20}) holds when each $F_i$ is
a star with a center at $V-U$
and $E(F_i)\subseteq E_G(U)$.

Assume that each $F_i$ is a star with a center at $V-U$
and $E(F_i)\subseteq E_G(U)$.
Then $G-U$ is an empty graph, 
$V=U\cup N_G(U)$ and $F$ is the bipartite graph 
$G[E_G(U)]$.

Let $E'=\{e_1,e_2,\cdots,e_t\}$, 
where $e_i$ is an edge in $F_i$. 
Applying Lemma~\ref{le2-01}~\ref{le2-1-n4}
and~\ref{le2-1-n3}  repeatedly, we have 
\begin{equation}\relabel{eq-pro2-1}
\tau_G(W)=\tau_{G/E'}(E(F/E'))
=\tau_{G[U]}(E(F/E')),
\end{equation}
where $F/E'$ is considered as 
a subforest of $G[U]$ 
whose components's vertex sets are 
$U\cap V(F_i)$ for $i=1,2,\cdots,t$.
Note that $F_0=G\spann{E(F/E')}$ 
is a spanning forest of $G[U]$
with 
$|U|-n_1-n_2-\cdots-n_t+t$ components 
with the following orders:
\begin{equation}\relabel{eq-pro2-2}
n_1, \cdots, n_t,
\underbrace{1,1,\cdots,1}_{|U|-t\ \mbox{numbers}}.
\end{equation} 
By Theorem~\ref{th1-1}, we have 
\begin{equation}\relabel{eq-pro2-3}
\tau_{G[U]}(F_0)=|U|^{|U|-2+t-n_1-\cdots-n_t} 
\prod_{i=1}^t n_i.
\end{equation}
Hence Claim~\thecountclaim\ holds.

\inclaim (\ref{eq3-20}) holds when each $F_i$ is
a star with a center at $V-U$.

Assume that each $F_i$ is a star with a center at $V-U$,
as shown in Figure~\ref{f13}.
\begin{figure}[h!]
 \centering
\unitlength 3mm 
\linethickness{0.4pt}
\ifx\plotpoint\undefined\newsavebox{\plotpoint}\fi 
\begin{picture}(22.937,14.982)(0,0)
\put(6.055,4.95){\circle*{.94}}
\put(17.103,4.95){\circle*{.94}}
\put(9.281,4.95){\circle*{.94}}
\multiput(7.071,10.828)(-.011168421,-.060484211){95}{\line(0,-1){.060484211}}
\multiput(18.12,10.828)(-.011168421,-.060484211){95}{\line(0,-1){.060484211}}
\put(9.281,4.906){\line(1,2){2.784}}
\multiput(6.055,4.95)(-.011199095,.02579638){221}{\line(0,1){.02579638}}
\multiput(17.103,4.95)(-.011199095,.02579638){221}{\line(0,1){.02579638}}
\multiput(6.01,4.95)(-.011204225,.039528169){142}{\line(0,1){.039528169}}
\multiput(17.059,4.95)(-.011204225,.039528169){142}{\line(0,1){.039528169}}
\multiput(9.237,4.95)(-.011083,.46775){12}{\line(0,1){.46775}}
\put(5.789,10.076){\makebox(0,0)[cc]{$\cdots$}}
\put(16.838,10.076){\makebox(0,0)[cc]{$\cdots$}}
\put(11.049,10.209){\makebox(0,0)[cc]{$\cdots$}}
\put(13.656,12.33){\makebox(0,0)[cc]{$U$}}
\put(12.949,4.773){$\cdots$}
\put(11.468,11.712){\oval(22.937,6.541)[]}
\put(3.978,5.966){\makebox(0,0)[cc]{\small {$F_1$}}}
\put(8.176,5.922){\makebox(0,0)[cc]{\small{$F_2$}}}
\put(18.606,6.01){\makebox(0,0)[cc]{\small{$F_t$}}}
\multiput(6.099,4.906)(.01105556,-.07244444){36}{\line(0,-1){.07244444}}
\multiput(9.237,4.994)(.011161616,-.027676768){99}{\line(0,-1){.027676768}}
\multiput(15.998,2.607)(.011161616,.024111111){99}{\line(0,1){.024111111}}
\multiput(17.103,4.994)(.011189873,-.017898734){158}{\line(0,-1){.017898734}}
\put(5.215,2.342){\makebox(0,0)[cc]{$\cdots$}}
\put(9.369,2.166){\makebox(0,0)[cc]{$\cdots$}}
\put(17.413,2.342){\makebox(0,0)[cc]{$\cdots$}}
\multiput(5.966,4.95)(-.011199095,-.011800905){221}{\line(0,-1){.011800905}}
\multiput(9.281,4.861)(-.011237288,-.02209322){118}{\line(0,-1){.02209322}}
\multiput(9.276,4.914)(.01113486,.09575978){59}{\line(0,1){.09575978}}
\end{picture}


\caption{Each $F_i$ is a star with a center in $V-U$}
\relabel{f13}
\end{figure}
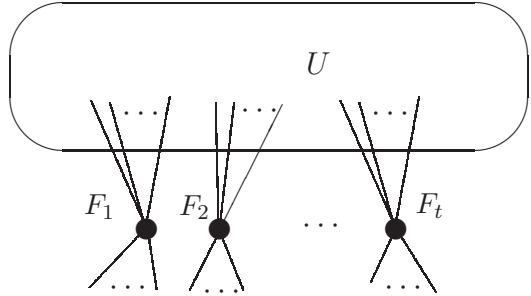
If $G-U$ is an independent set of $G$, 
then each $F_i$ is a star with a center in $V-U$ 
and $E(F_i)\subseteq E_G(U)$, 
implying that the claim holds by Claim 1. 

Now assume that $V-U$ is not independent in $G$,
i.e., $E_0=E(G-U)\ne \emptyset$.
Since each component $F_i$ is a star with a center in $V-U$, 
each edge $e\in E_0$ is incident with two vertices in $V-U$
one of which is an end-vertex.
Thus $U$ is still a clique of $G/E_0$ and 
$F/E_0$ is a forest 
with $t$ components $F'_1, F'_2, \cdots, F'_t$
each of which is a star with 
a center in $V(G/E_0)-U$ and 
each edge in $E(F/E_0)$ is incident with some vertex in $U$,
where $F'_i=F_i/(E_0\cap E(F_i))$. 
By Claim 1, the result holds for $G/E_0$, i.e., 
\begin{equation}\relabel{eq3-20-2}
\tau_{G/E_0}(E(F/E_0))=
|U|^{|U|-2+t-n'_1-n'_2-\cdots-n'_t}
\prod_{i=1}^t n'_i,
\end{equation}
where $n'_i=|U\cap V(F'_i)|$.
Clearly $n'_i=|U\cap V(F'_i)|=|U\cap V(F_i)|=n_i$.
Applying Lemma~\ref{le2-01}~\ref{le2-1-n4} repeatedly,
we have 
$\tau_G(E(F))=\tau_{G/E_0}(E(F/E_0))$.
Thus Claim 2 holds.

\inclaim (\ref{eq3-20}) holds whenever $F$ is a forest.

Let $W=E(F)=E(G-G[U])$
and $G'=G\star W-W$. 
Observe that $U$ is a clique of $G'$ 
and $G'[E(G')-E(G'[U])]$ is a forest 
with $t$ components $F'_1,\cdots, F_t'$ 
each of which is a star with a center in $V(G')-U$
such that $V(F'_i)\cap U=V(F_i)\cap U$ holds for all 
$i=1,2,\cdots,t$.
By Claim 2, the result holds 
for $G'$, i.e., 
\begin{equation}\relabel{eq3-20-3}
\tau_{G'}(W')=
|U|^{|U|-2+t-n_1-n_2-\cdots-n_t}
\prod_{i=1}^t n_i,
\end{equation}
where $W'=E(G\star W)-W$.
By Lemma~\ref{le2-11}, we have 
$\tau_G(W)=\tau_{G\star W-W}(W')$.

Thus Claim 3 holds and the result is proved. 
\proofend

\noindent {\bf Remark}: 
For any forest $M$ in $K_n$ with components $M_1,\cdots,M_t$, 
let $G=K_n\star M$. 
By Lemma~\ref{le2-11}, 
$\tau_{K_n}(M)=\tau_{G}(W)$, 
where $W=E(G)-E(K_n)$.
Note that each component $F_i$ of $G[W]$ is 
a star with a center in $V(G)-U$ 
and $E(F_i)\subseteq E_G(U)$, where $U=V(K_n)$
and $V(F_i)\cap U=V(M_i)$.
Theorem~\ref{th1-1} corresponds to 
Proposition~\ref{pro3-1} for the case that 
each component $F_i$ of $F$ is a star with a center 
in $V-U$ and $E(F_i)\subseteq E_G(U)$.

\subsection{
Relation between $\tau_G(M\cup N)$ 
and $\tau_{G\bullet U}(N)$
\relabel{sec3-2}
}

In this subsection, we assume that 
$U$ is a clique of $G$ and each component of $G[E_G(U)]$ 
is a star with a center in $V-U$.
Note that  each component of $G[E_G(U)]$ is a star 
with a center in $V-U$ 
\iff each vertex in $U$ is incident 
with at most one edge in $E_G(U)$.

Let $u$ be the new vertex in $G\bullet U$
created after contracting all edges in $E(G[U])$.
So the vertex set of $G\bullet U$ is $(V-U)\cup \{u\}$.
Note that $G\bullet U$ may have parallel edges incident with
vertex $u$, as all edges in $E_G(U)=E_G(U, V-U)$
are the edges in $G\bullet U$ incident with $u$.
For each $v\in V-U$, 
the number of parallel edges in $G\bullet U$ 
joining $u$ and $v$ is equal to $|N_G(v)\cap U|$.

Now we are going to establish 
the main result in this section.

\begin{theo}\relabel{th3-1} 
Let $M=E_G(U)$ and $N\subseteq E(G-U)$.
If each component of $G[M]$ is a star with a center in $V-U$,  
then 
\begin{equation}\relabel{eq3-2}
\tau_G(M\cup N)=
|U|^{|U|-2-|M|}
\sum_{T\in \setst_{G\bullet U}(N)}
|U|^{|E_T(u)|}.
\end{equation}
\end{theo}

\proof By the given condition on $M$, 
$G[M\cup N]$ contains cycles \iff 
$G\bullet U[N]$ contains cycles, 
implying that (\ref{eq3-2}) holds whenever 
$G[M\cup N]$ contains cycles.
Thus, it suffices to consider the case that 
$G[M\cup N]$ is a forest.
We will prove (\ref{eq3-2}) 
by completing the following claims. 

\inclaim 
(\ref{eq3-2}) holds if
$G-E(G[U])$ is a forest and $N=E(G-U)$.

\begin{figure}[htbp]
 \centering
\unitlength 3mm 
\linethickness{0.4pt}
\ifx\plotpoint\undefined\newsavebox{\plotpoint}\fi 
\begin{picture}(23.203,12.64)(0,0)
\put(5.436,1.79){\oval(6.01,3.933)[]}
\put(16.485,1.79){\oval(6.01,3.933)[]}
\put(4.154,2.652){\circle*{.94}}
\put(15.203,2.652){\circle*{.94}}
\put(7.38,2.652){\circle*{.94}}
\put(18.429,2.652){\circle*{.94}}
\multiput(5.171,8.53)(-.011168421,-.060484211){95}{\line(0,-1){.060484211}}
\multiput(16.219,8.53)(-.011157895,-.060484211){95}{\line(0,-1){.060484211}}
\put(7.38,2.607){\line(1,2){2.784}}
\put(18.429,2.607){\line(1,2){2.784}}
\put(5.745,2.475){\makebox(0,0)[cc]{$\cdots$}}
\put(16.794,2.475){\makebox(0,0)[cc]{$\cdots$}}
\multiput(4.154,2.652)(-.011199095,.02579638){221}{\line(0,1){.02579638}}
\multiput(15.203,2.652)(-.011199095,.02579638){221}{\line(0,1){.02579638}}
\multiput(4.11,2.652)(-.011204225,.039521127){142}{\line(0,1){.039521127}}
\multiput(15.159,2.652)(-.011204225,.039521127){142}{\line(0,1){.039521127}}
\multiput(7.336,2.652)(-.011,.467667){12}{\line(0,1){.467667}}
\multiput(18.385,2.652)(-.011083,.467667){12}{\line(0,1){.467667}}
\put(3.889,7.778){\makebox(0,0)[cc]{$\cdots$}}
\put(14.938,7.778){\makebox(0,0)[cc]{$\cdots$}}
\put(8.971,7.866){\makebox(0,0)[cc]{$\cdots$}}
\put(20.02,7.866){\makebox(0,0)[cc]{$\cdots$}}
\multiput(7.425,2.652)(.0110313,.175375){32}{\line(0,1){.175375}}
\multiput(18.473,2.652)(.0110625,.175375){32}{\line(0,1){.175375}}
\put(11.756,10.032){\makebox(0,0)[cc]{$U$}}
\put(10.386,1.635){$\cdots$}
\put(11.734,9.369){\oval(22.937,6.541)[]}
\put(.972,3.094){\makebox(0,0)[cc]{$F_1$}}
\put(20.86,3.005){\makebox(0,0)[cc]{$F_t$}}
\end{picture}


\caption{$G[M\cup N]$ is forest with $t$ components 
$F_1,\cdots, F_t$}
\relabel{f12}
\end{figure}
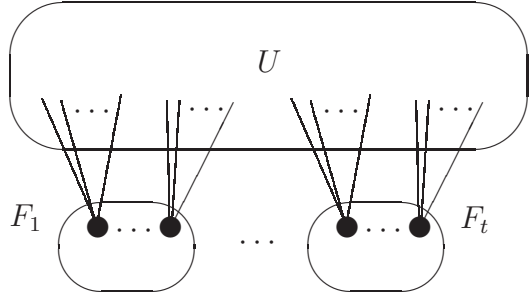

Assume that $N=E(G-U)$
and $G[M\cup N]$ is a forest
with components $F_1, F_2, \cdots, F_t$,
as shown in Figure~\ref{f12}.
By Proposition~\ref{pro3-1}, 
\begin{equation}\relabel{eq3-6}
\tau_G(M\cup N)=
|U|^{|U|-2+t-n_1-\cdots-n_t}
\prod_{i=1}^t n_i,
\end{equation}
where 
$n_i=|V(F_i)\cap U|$ for $i=1,2,\cdots,t$, 
implying that (\ref{eq3-2}) holds if and only if the following equality holds:
\begin{equation}\relabel{eq3-7}
|U|^{t-n_1-\cdots-n_t}\prod_{i=1}^t n_i
=U^{-|M|}\sum_{T\in \setst_{G\bullet U}(N)}
|U|^{|E_T(u)|}.
\end{equation}
By the given condition, each vertex in $U$ 
is incident with at most one edge in $M$.
Since $n_1+n_2+\cdots+n_t$ is the number of 
vertices in $U$ which are incident with edges in $M$,
we have $n_1+n_2+\cdots+n_t=|M|$.
For any $T\in \setst_{G\bullet U}(N)$, we have 
$|E(T)\cap E_G(U, V(F_i)-U)|=1$ for all $i=1,2,\cdots,t$,
implying that $|E_T(u)|=t$.

It remains to show that 
$\tau_{G\bullet U}(N)=\prod_{i=1}^t n_i$.
Let $T\in \setst_{G\bullet U}(N)$.
Observe that $T-u$ is actually the graph $G-U$, 
which consists of $t$ 
components $F_i-V(F_i)\cap U$
for $i=1,2,\cdots,t$.
Also note that $T$ contains exactly $t$ edges 
$e_1,e_2,\cdots,e_t$, where each $e_i$
with $u$ and some vertex in $F_i-V(F_i)\cap U$
for $i=1,2,\cdots,t$.
Observe that each $e_i$ can be any one of the edges 
in the set $M\cap E(F_i)$ whose size is exactly 
$|V(F_i)\cap M|=n_i$.
Hence 
\begin{equation}\relabel{eq3-7-1}
\tau_{G\bullet U}(N)
=\prod_{i=1}^t |V(F_i)\cap U|
=\prod_{i=1}^t n_i.
\end{equation}
Thus (\ref{eq3-7}) holds and 
Claim~\thecountclaim\ follows.

\inclaim 
(\ref{eq3-2}) holds for any $N\subseteq E(G-U)$ 
such that $G[M\cup N]$ is a forest.

For any $T\in \setst_G(M\cup N)$, $T-U$ is a forest 
with $N\subseteq E(T-U)$.
Let $E_0=E(G(U))$ and let 
$\setn$ be the family of those 
subsets $N'$ of $E(G-U)$ with 
$N\subseteq N'$ such that 
$G\spann{N'}$ is 
a forest and 
$G\spann {M\cup N'\cup E_0}$ is connected. 

Clearly,  
$\setst_{G\spann {M\cup N_1\cup E_0}}(M\cup N_1)$
and $\setst_{G\spann {M\cup N_2\cup E_0}}(M\cup N_2)$
are disjoint for any pair of distinct members $N_1,N_2\in \setn$,
and 
\begin{equation}\relabel{eq3-3}
\setst_G(M\cup N)=\bigcup_{N'\in \setn} 
\setst_{G\spann {M\cup N'\cup E_0}}(M\cup N').
\end{equation}
Similarly, for any pair of distinct members 
$N_1,N_2\in \setn$,
$\setst_{G\bullet U\spann {M\cup N_1}}(N_1)$
and $\setst_{G\bullet U\spann {M\cup N_2}}(N_2)$ 
are disjoint, and 
\begin{equation}\relabel{eq3-4}
\setst_{G\bullet U}(N)=\bigcup_{N'\in \setn} 
\setst_{G\bullet U\spann {M\cup N'}}(N').
\end{equation}
By Claim 1, 
the following identity holds for any $N'\in \setn$: 
\begin{equation}\relabel{eq3-5}
\tau_{G\spann {M\cup N'\cup E_0}}(M\cup N')
=|U|^{|U|-2-|M|}
\sum_{T\in \setst_{G\bullet U\spann {M\cup N'}}(N')}
|U|^{|E_T(u)|}.
\end{equation}
Thus Claim~\thecountclaim\ follows from (\ref{eq3-3}), (\ref{eq3-4})
and (\ref{eq3-5}).  
\proofend

\resection{
When $V_1,V_2,\cdots,V_k$ are 
disjoint cliques of $G$
\relabel{sec4}}

In this section, we always assume that 
$G=(V,E)$ is a connected and loopless multigraph,
where $V$ is partitioned into non-empty subsets 
$V_0,V_1,\cdots, V_k$ 
satisfying the following conditions:
\begin{enumerate}
\item\relabel{sec5-con0} $V_i$ is a clique for all $i=1,2,\cdots,k$;
\item\relabel{sec5-con1}  
$E_G(V_i,V_j)=\emptyset$ for each pair $i,j$ with 
$1\le i<j\le k$; 
\item\relabel{sec5-con2} 
each component of $G[M]$ is a star with a center in $V_0$, 
where 
$M=\bigcup\limits_{1\le i\le k}E_G(V_0,V_i)$.
\end{enumerate}

The structure of $G$ under conditions (i), (ii) and (iii) above  
is as shown in Figure~\ref{f4}(a).
Note that condition~\ref{sec5-con2} above is equivalent to 
that each vertex in $V_i$ is incident with at most one 
edge in $M$ for all $i=1,2,\cdots,k$.
All parallel edges of $G$ must be in the 
subgraph $G[V_0]$.

\begin{figure}[http]
 \centering

\unitlength .6mm 
\linethickness{0.4pt}
\ifx\plotpoint\undefined\newsavebox{\plotpoint}\fi 

  
(a) $G$ \hspace{6.5 cm} (b) $G\bullet U$

\caption{$G-V_0$ consists of $k$ cliques and 
$G[M]$ consists of stars with centers in $V_0$}
\relabel{f4}
\end{figure}

Let $M_i=E_G(V_i,V_0)$ for all $i=1,2,\cdots,k$. 
Then $M=E_G(V_0)=\bigcup\limits_{1\le i\le k}M_i$.
In this section,   
our main purpose is to apply
Theorem~\ref{th3-1} to find an 
expression for $\tau_G(M\cup N)$
for any $N\subseteq E(G[V_0])$.
Applying this result, we are able to get 
an expression of $\tau_G(R\cup N)$
for any $R\subseteq M$.

\subsection{$\tau_G(M\cup N)$ for $N\subseteq E(G[V_0])$
\relabel{sec4-1}}

Let $U=V_1\cup \cdots \cup V_k$. 
Recall that 
$G\bullet U$ is defined to be the graph $G/E(G[U])$.
As $G[U]$ has $k$ components $G[V_1], G[V_2],
\cdots,G[V_k]$,  
$G\bullet U$ can be obtained from $G$ 
by removing all edges in $G[U]$ and 
identifying all vertices in each $V_i$ 
as one vertex, denoted by $v_i$, for $i=1,2,\cdots,k$,
as shown in  Figure~\ref{f4} (b).
Thus $V(G\bullet U)=V_0\cup \{v_1,v_2,\cdots,v_k\}$
and $E(G\bullet U)
=E(G)-\bigcup_{1\le i\le k}E(G[V_i])=M\cup E(G[V_0])$.

\begin{theo}\relabel{th4-2}  
For any $N\subseteq E(G[V_0])$, 
\begin{equation}\relabel{eq4-2}
\tau_G(M\cup N)=
\prod_{i=1}^k |V_i|^{|V_i|-2-|M_i|}
\sum_{T\in \setst_{G\bullet U}(N)}
\prod_{i=1}^k 
|V_i|^{|E_T(v_i)|},
\end{equation} 
where $M_i=E_G(V_i,V_0)$ for $i=1,2,\cdots,k$
and $M=M_1\cup \cdots \cup M_k$.
\end{theo}

\proof
If $k=1$, the result follows directly from 
Theorem~\ref{th3-1}.
Assume that the result holds for $k<n$, where $n\ge 2$.
Now consider the case that $k=n$.

By Theorem~\ref{th3-1},
\begin{eqnarray}
\frac{\tau_G(M\cup N)}{|V_k|^{|V_k|-2-|M_k|}}
&=&
\sum_{T\in \setst_{G\bullet V_k}((M-M_k)\cup N)}|V_k|^{|E_T(v_k)|}\nonumber \\
&=&
\sum_{B_k\subseteq M_k}
\sum_{T\in \setst_{G\bullet V_k}((M-M_k)\cup N)
\atop E_T(v_k)=M_k-B_k}
|V_k|^{|M_k-B_k|}\nonumber \\
&=&
\sum_{B_k\subseteq M_k}|V_k|^{|M_k-B_k|}
\tau_{(G-B_k)\bullet V_k} ((M-B_k)\cup N),
\relabel{eq4-001}
\end{eqnarray}
where the last equality follows from the fact that 
$T\in \setst_{G\bullet V_k}((M-M_k)\cup N)$ 
with $E_T(v_k)=M_k-B_k$ 
\iff $T\in \setst_{(G-B_k)\bullet V_k} ((M-B_k)\cup N)$,
as $M-B_k=(M-M_k)\cup (M_k-B_k)$.

For any $B_k\subseteq M_k$,
we have 
$M-B_k=M_1\cup \cdots \cup M_{k-1}\cup (M_k-B_k)$.
By the inductive assumption,
\begin{equation}
\frac{\tau_{(G-B_k)\bullet V_k} ((M-B_k)\cup N)}
{\prod_{i=1}^{k-1} |V_i|^{|V_i|-2-|M_i|}}
=
\sum_{T\in \setst_{(G-B_k)\bullet U}(N\cup (M_k-B_k))}
\prod_{i=1}^{k-1} 
|V_i|^{|E_T(v_i)|}.\relabel{eq4-002}
\end{equation} 
By (\ref{eq4-001}) and (\ref{eq4-002}), 
we have 
\begin{eqnarray}
\frac{\tau_G(M\cup N)}
{\prod_{i=1}^{k} |V_i|^{|V_i|-2-|M_i|}}
&=&\sum_{B_k\subseteq M_k}
\sum_{T\in \setst_{(G-B_k)\bullet U}(N\cup (M_k-B_k))}
|V_k|^{|M_k-B_k|}
\prod_{i=1}^{k-1} 
|V_i|^{|E_T(v_i)|}\nonumber \\
&=&
\sum_{B_k\subseteq M_k}
\sum_{T\in \setst_{(G-B_k)\bullet U}(N)
\atop E_T(v_k)=M_k-B_k}
|V_k|^{|E_T(v_k)|}
\prod_{i=1}^{k-1} 
|V_i|^{|E_T(v_i)|}\nonumber \\
&=&
\sum_{T\in \setst_{G\bullet U}(N)}
\prod_{i=1}^{k} 
|V_i|^{|E_T(v_i)|}. \relabel{eq4-003}
\end{eqnarray}
\proofend

Now we give a proof of Theorem~\ref{th1-2} 
by applying Theorem~\ref{th4-2} directly.

\noindent {\it Proof} of Theorem~\ref{th1-2}.
Assume that $G=(V,E)$ is a graph satisfying 
the conditions assumed in the beginning of this section,
$V_0$ is an independent set of $G$ and 
each component of $G[M]$ is a star of size $2$.

Let $V_0=\{w_1,w_2,\cdots,w_r\}$. 
Then $G[M]$ consists of exactly $r$ components
$S_1, S_2,\cdots, S_r$ 
which are stars of size $2$ 
with centers $w_1,w_2,\cdots, w_r$
respectively.
Let  $e_{i,1}$ and $e_{i,2}$ denote the two edges
in $S_i$.

Clearly, 
$\{e_{i,1},e_{i,2}\}\cap E(T)\ne \emptyset$
holds for any $1\le i\le r$ and 
any $T\in \setst_{G\bullet U}$.
For any $I\subseteq \{1,2,\cdots,r\}$,
let $\setst_{G\bullet U}^{I}$ be the set of 
members $T\in \setst_{G\bullet U}$ 
such that 
$$
\{1\le i\le r: \{e_{i,1},e_{i,2}\}\subseteq E(T)\}=I.
$$ 
Observe that $\setst_{G\bullet U}^{I}\ne \emptyset$
\iff the edge set $\{e_{i,1},e_{i,2}:i\in I\}$ 
induces a spanning tree of the subgraph 
$(G\bullet U)-\{w_j: 1\le j\le r, j\notin I\}$.
When $\setst_{G\bullet U}^{I}\ne \emptyset$,
there are exactly $2^{r-|I|}$ members in
$\setst_{G\bullet U}^{I}$
and for each $T\in \setst_{G\bullet U}^{I}$, 
$|E(T)\cap \{e_{i,1},e_{i,2}\}|=1$ holds 
for all $i\in \{1,2,\cdots,r\}-I$.

Let $G'$ be the graph obtained from $G$ by 
contracting exactly one edge in each $S_i$ 
for all $i=1,2,\cdots,r$,
and let $M'=M\cap E(G')$.
By Lemma~\ref{le2-01}~\ref{le2-1-n4},
$\tau_{G'}(M')=\tau_{G}(M)$.

For any $e\in M$, let $l(e)=i$ 
such that $e$ is incident with a vertex in $V_i$.
By Theorem~\ref{th4-2}, 
\begin{eqnarray}
\tau_{G}(M)
&=&\prod_{i=1}^k |V_i|^{|V_i|-2-|M_i|}
\sum_{T\in \setst_{G\bullet U}}\prod_{i=1}^k |V_i|^{|E_T(v_i)|}
\nonumber \\
&=&\prod_{i=1}^k |V_i|^{|V_i|-2}
\sum_{I\subseteq \{1,2,\cdots,r\}}
\sum_{T\in \setst_{G\bullet U}^I}
\prod_{i=1}^k |V_i|^{-|M_i-E_T(v_i)|}
\nonumber \\
&=&\prod_{i=1}^k |V_i|^{|V_i|-2}
\sum_{I\subseteq \{1,2,\cdots,r\}}
\sum_{T\in \setst_{G\bullet U}^I}
\prod_{e\in M-E(T)}|V_{l(e)}|^{-1}
\nonumber \\
&=&\prod_{i=1}^k |V_i|^{|V_i|-2}
\sum_{I\subseteq \{1,2,\cdots,r\}
\atop \setst_{G\bullet U}^I\ne \emptyset}
\prod_{i\in \{1,2,\cdots,r\}-I}
(|V_{l(e_{i,1})}|^{-1}+|V_{l(e_{i,2})}|^{-1})
\nonumber \\
&=&\prod_{i=1}^k |V_i|^{|V_i|-2}
\sum_{T'\in \setst_{G^*}}
\prod_{e\in M'-E(T')}
(|V_{a(e)}|^{-1}+|V_{b(e)}|^{-1}),
\relabel{eq4-5}
\end{eqnarray}
where $G^*$ is the graph 
obtained from $G'$ by identifying 
all vertices in each $V_i$ as a vertex, denoted by $v_i$,
and removing all loops, 
and $a(e)$ and $b(e)$ 
are numbers in $\{1,2,\cdots,k\}$ such that 
$v_{a(e)}$ and $v_{b(e)}$ are the two ends 
of $e$ in $G^*$
which correspond to $V_{a(e)}$ and $V_{b(e)}$.

Since $|\setst_{G'}(M')|=\tau_{G}(M)$,
Theorem~\ref{th1-2} is proven by (\ref{eq4-5}). 
\proofend

\subsection{$\tau_G(R\cap N)$ for $R\subseteq M$
and $N\subseteq E(G-U)$
\relabel{sec4-2}}

In this subsection, we will find an 
expression for $\tau_G(R\cup N)$
for any $R\subseteq E_G(U)$ and $N\subseteq E(G-U)$.

\begin{theo}\relabel{th5-3}
For any $R\subseteq M$ and $N\subseteq E(G-U)$,
\begin{eqnarray}\relabel{eq5-6}
\tau_G(R\cup N)
=
\prod_{i=1}^k |V_i|^{|V_i|-2}
\sum_{T\in \setst_{G\bullet U}(N)}
\prod_{i=1}^k |V_i|^{-|M_i- E(T)|}
(1+|V_i|)^{|(M_i-R)- E(T)|}.
\end{eqnarray}

\end{theo}

\proof Let $R$ be a fixed subset of $M$.
Note that $T\in \setst_G(R\cup N)$ 
\iff $T\in \setst_{G-B}((M-B)\cup N)$ for 
some $B$ with $B\subseteq M-R$.
Obviously, for distinct subsets $B_1, B_2$ of $M-R$, 
$\setst_{G-B_1}((M-B_1)\cup N)$ and 
$\setst_{G-B_2}((M-B_2)\cup N)$ are disjoint. 
Thus,
\begin{equation}\relabel{eq5-6-1}
\tau_G(R\cup N) =
\sum_{B\subseteq M-R}
\tau_{G-B}((M-B)\cup N).
\end{equation}
Then, by Theorem~\ref{th4-2}, 
\begin{eqnarray}
\tau_G(R\cup N)   
&=&
\sum_{B\subseteq M-R}
\sum_{T\in \setst_{(G\bullet U)-B}(N)}
\prod_{i=1}^k |V_i|^{|V_i|-2-|M_i-B|+
|(M_i-B)\cap E_T(v_i)|}
\nonumber \\ &=&
\prod_{i=1}^k |V_i|^{|V_i|-2}
\sum_{B\subseteq M-R}
\sum_{T\in \setst_{(G\bullet U)-B}(N)}
\prod_{i=1}^k |V_i|^{-|M_i-B|+
|(M_i-B)\cap E_T(v_i)|}
\nonumber \\ &=&
\prod_{i=1}^k |V_i|^{|V_i|-2}
\sum_{B\subseteq M-R}
\sum_{R'\subseteq M-B}
\sum_{T\in \setst_{(G\bullet U)-B}(N)
\atop E(T)\cap M=R'}
\prod_{i=1}^k |V_i|^{-|M_i-B|+
|R'\cap M_i|}
\nonumber \\ &=&
\prod_{i=1}^k |V_i|^{|V_i|-2}
\sum_{R'\subseteq M}
\prod_{i=1}^k |V_i|^{|R'\cap M_i|}
\Phi(R')
\relabel{eq5-19}
\end{eqnarray}
where 
\begin{eqnarray}
\Phi(R') &=&
\sum_{B\subseteq M-(R'\cup R)}
\sum_{T\in \setst_{(G\bullet U)-B}(N)
\atop E(T)\cap M=R'}
\prod_{i=1}^k |V_i|^{-|M_i-B|}
\nonumber \\
&=&
\prod_{i=1}^k |V_i|^{-|R\cap M_i|}
\sum_{B\subseteq M-(R'\cup R)}
\sum_{T\in \setst_{(G\bullet U)-B}(N)
\atop E(T)\cap M=R'}
\prod_{i=1}^k |V_i|^{-|(M_i-B)-R|}
\nonumber \\
&=&
\prod_{i=1}^k |V_i|^{-|R\cap M_i|}
\sum_{T\in \setst_{G\bullet U}(N)
\atop E(T)\cap M=R'}
\sum_{B\subseteq M-(R'\cup R)}
\prod_{i=1}^k |V_i|^{-|(M_i-R)-B|}
\nonumber \\
&=&
\prod_{i=1}^k |V_i|^{-|(R\cup R')\cap M_i|}
\sum_{T\in \setst_{G\bullet U}(N)
\atop E(T)\cap M=R'}
\sum_{B\subseteq M-(R'\cup R)}
\prod_{i=1}^k |V_i|^{-|M_i-(R\cup R')-B|}
\nonumber \\
&=&
\prod_{i=1}^k |V_i|^{-|(R\cup R')\cap M_i|}
\sum_{T\in \setst_{G\bullet U}(N)
\atop E(T)\cap M=R'}
\prod_{i=1}^k (1+|V_i|^{-1})^{|M_i-(R\cup R')|}.
\relabel{eq5-19-1}
\end{eqnarray}
By (\ref{eq5-19-1}) and (\ref{eq5-19}), we have 
\begin{eqnarray}
\frac{\tau_G(R\cup N)}
{\prod_{i=1}^k |V_i|^{|V_i|-2}}
&=&
\sum_{R'\subseteq M}
\sum_{T\in \setst_{G\bullet U}(N)\atop E(T)\cap M=R'}
\prod_{i=1}^k |V_i|^{-|(R-R')\cap M_i|}
(1+|V_i|^{-1})^{|M_i-(R\cup R')|}
\nonumber \\ 
&=& 
\sum_{T\in \setst_{G\bullet U}(N)}
\prod_{i=1}^k |V_i|^{-|(M_i\cap R)-E(T)|}
(1+|V_i|^{-1})^{|(M_i-R)-E(T))|}.
\relabel{eq5-18}
\end{eqnarray}
Thus we can verify that the result holds.
\proofend

When $R=\emptyset$, a direct application of 
Theorem~\ref{th5-3} gives an expression for 
$\tau_{G}(N)$.

\begin{cor}\relabel{co5-3-1}
For any $N\subseteq E(G-U)$,
\begin{eqnarray}\relabel{co5-3-1-ea1}
\tau_G(N)
=
\prod_{i=1}^k |V_i|^{|V_i|-2}
 \sum_{T\in \setst_{G\bullet U}(N)}
\prod_{i=1}^k 
(1+1/|V_i|)^{|M_i-E(T)|}.
\end{eqnarray}
\end{cor}

\subsection{
Simplify the expression of (\ref{eq5-6}) 
\relabel{sec4-3}
}

Note that (\ref{eq5-6}) can be changed to
\begin{eqnarray}\relabel{eq5-6-0}
\tau_G(R\cup N)
&=&
\prod_{i=1}^k |V_i|^{|V_i|-2-|M_i|}
(1+|V_i|)^{|M_i-R|}
\nonumber \\
& &
\times \sum_{T\in \setst_{G\bullet U}(N)}
\prod_{i=1}^k |V_i|^{|M_i\cap R\cap E(T)|}
\left (\frac{|V_i|}{1+|V_i|}\right )^{|(M_i-R)\cap E(T)|}.
\end{eqnarray}

For any $R\subseteq M$, 
let $\omega$ be the mapping from 
$E(G\bullet U)$ (i.e., $M\cup E(G-U)$)
to $\N=\{1,2,3,\cdots\}$ defined below:
\begin{equation}\relabel{eq3-3-1}
\omega_R(e)=
\left \{
\begin{array}{ll}
|V_i|, \qquad &e\in M_i\cap R;\\
|V_i|/(1+|V_i|), \qquad &e\in M_i-R;\\
1,  &\mbox{otherwise}.
\end{array}
\right.
\end{equation}
Then (\ref{eq5-6}) can be expressed as
\begin{equation}
\tau_G(R\cup N)
=\prod_{i=1}^k |V_i|^{|V_i|-2-|M_i|}
(1+|V_i|)^{|M_i-R|}
\sum_{T\in \setst_{G\bullet U}(N)}
\prod_{e\in E(T)}\omega_R(e). \relabel{eq3-3-2}
\end{equation}
Let $w_1,w_2,\cdots,w_r$ be the vertices in the set 
$N_G(U)-U$. 
These vertices are actually centers 
of the components of $G[M]$, 
as each component of $G[M]$ is a star.

As there may be more than one edge in 
$E_{G\bullet U}(v_i,w_j)\cap R$ 
or 
$E_{G\bullet U}(v_i,w_j)-R$  for 
$1\le i\le k$ and $1\le j\le r$,
(\ref{eq3-3-2}) can be further simplified. 

Given $R\subseteq M$, 
let $G\circ_R  U$ denote the graph 
obtained from $G\bullet U$ by removing 
$|E_{G\bullet U}(v_i,w_j)\cap R|-1$ edges 
in the set $E_{G\bullet U}(v_i,w_j)\cap R$ 
whenever $|E_{G\bullet U}(v_i,w_j)\cap R|\ge 2$
and  $|E_{G\bullet U}(v_i,w_j)-R|-1$ edges 
in the set $E_{G\bullet U}(v_i,w_j)-R$
whenever $|E_{G\bullet U}(v_i,w_j)-R|\ge 2$
for each pair $i,j:1\le i\le k$ and $1\le j\le r$.
Thus, in the graph $G\circ_R  U$, 
there are at most two parallel 
edges joining each pair of vertices $v_i$ and $w_j$.
If this case happens, then exactly one 
of the two edges joining $v_i$ and $w_j$ is contained in $R$. 
 
Let $\omega'_R$ be the mapping from 
$E(G\circ_R  U)$ 
to $\N=\{1,2,3,\cdots\}$ defined below:
\begin{equation}\relabel{eq3-3-3}
\omega_R'(e)=
\left \{
\begin{array}{ll}
|V_i|\cdot |E_{G\bullet U}(v_i,w_j)\cap R|, 
&e\in R\ \mbox{ and }  
e\mbox{ joins }v_i \mbox{ and }w_j;\\
\frac{|V_i|}{1+|V_i|} \cdot |E_{G\bullet U}(v_i,w_j)-R|, 
&e\notin R\ \mbox{ and }  
e\mbox{ joins }v_i \mbox{ and }w_j;\\
1,  &\mbox{otherwise},
\end{array}
\right.
\end{equation}
where $1\le i\le k$ and $1\le j\le r$.
Then (\ref{eq3-3-2}) can be replaced by 
the following expression:
\begin{equation}\relabel{eq3-3-4}
\tau_G(R\cup N)
=\prod_{i=1}^k |V_i|^{|V_i|-2-|M_i|}
(1+|V_i|)^{|M_i-R|}
\sum_{T\in \setst_{G\circ_R U}(N)}
\prod_{e\in E(T)}\omega_R'(e).
\end{equation}

\resection{ 
When $E_1,E_2,\cdots,E_k$  
is a partition of $E$ 
such that each $G[E_i]$ is a complete graph
\relabel{sec5}
}

\subsection{$\tau_G=|\setst_{G\diamond \sets}|$ holds
for a graph $G\diamond \sets$
\relabel{sec5-01}}

Let $v$ be any vertex in $G$ and $E_0\subseteq E_G(v)$.
Let $G_{v\triangleleft   E_0}$ denote the graph obtained from 
$G-(E_G(v)-E_0)$ by adding a new vertex $v'$
and a new edge joining $v$ and $v'$ and 
finally changing the end $v$ of all edges in 
$E_G(v)-E_0$ to $v'$,
as shown in Figure~\ref{f8}.
Clearly $G\triangleleft   E_0/e'\cong G$,
where $e'=vv'$ is the only edge in 
$E(G\triangleleft E_0)-E(G)$.
By Lemma~\ref{le2-01}~\ref{le2-1-n4}, 
$\tau_G(W)=|\setst_{G_{v\triangleleft E_0}}(W\cup e')|$
holds.

\begin{figure}[h!]
 \centering

\unitlength 2mm 
\linethickness{0.4pt}
\ifx\plotpoint\undefined\newsavebox{\plotpoint}\fi 
\begin{picture}(44.001,13.23)(0,0)
\put(6.02,5.723){\circle*{1.27}}
\put(37.46,5.544){\circle*{1.27}}
\put(31.44,5.723){\circle*{1.27}}
\put(5.797,5.723){\line(1,0){.149}}
\put(31.365,5.575){\line(1,0){.149}}
\multiput(1.635,13.23)(.0168398438,-.0293242187){256}{\line(0,-1){.0293242187}}
\multiput(27.203,13.081)(.0168398438,-.0293242187){256}{\line(0,-1){.0293242187}}
\multiput(5.946,5.723)(-.0176595745,-.0168687943){282}{\line(-1,0){.0176595745}}
\multiput(31.514,5.574)(-.0176595745,-.0168652482){282}{\line(-1,0){.0176595745}}
\multiput(1.189,9.068)(.024902062,-.016860825){194}{\line(1,0){.024902062}}
\multiput(26.757,8.919)(.024907216,-.01685567){194}{\line(1,0){.024907216}}
\multiput(6.02,5.797)(.0168626866,.0213014925){335}{\line(0,1){.0213014925}}
\multiput(37.46,5.618)(.0168626866,.0213014925){335}{\line(0,1){.0213014925}}
\multiput(5.872,5.723)(.028163158,.016821053){190}{\line(1,0){.028163158}}
\multiput(37.312,5.544)(.028163158,.016821053){190}{\line(1,0){.028163158}}
\multiput(5.946,5.797)(.016848889,-.018826667){225}{\line(0,-1){.018826667}}
\multiput(37.386,5.618)(.016848889,-.018826667){225}{\line(0,-1){.018826667}}
\put(6.02,7.656){\makebox(0,0)[cc]{$v$}}
\put(37.46,7.477){\makebox(0,0)[cc]{$v'$}}
\put(3.865,11.075){\makebox(0,0)[cc]{$e_1$}}
\put(29.433,10.926){\makebox(0,0)[cc]{$e_1$}}
\put(2.899,8.622){\makebox(0,0)[cc]{$e_2$}}
\put(28.467,8.473){\makebox(0,0)[cc]{$e_2$}}
\put(1.414,2.409){\makebox(0,0)[cc]{$e_s$}}
\put(26.982,2.26){\makebox(0,0)[cc]{$e_s$}}
\put(12.561,11.075){\makebox(0,0)[cc]{$e'_1$}}
\put(44.001,10.896){\makebox(0,0)[cc]{$e'_1$}}
\put(12.041,7.433){\makebox(0,0)[cc]{$e'_2$}}
\put(43.481,7.254){\makebox(0,0)[cc]{$e'_2$}}
\put(11.521,2.453){\makebox(0,0)[cc]{$e'_t$}}
\put(42.961,2.274){\makebox(0,0)[cc]{$e'_t$}}
\put(9.811,5.574){\makebox(0,0)[cc]{$\vdots$}}
\put(41.251,5.395){\makebox(0,0)[cc]{$\vdots$}}
\put(1.932,5.5){\makebox(0,0)[cc]{$\vdots$}}
\put(27.501,5.351){\makebox(0,0)[cc]{$\vdots$}}
\put(31.589,7.433){\makebox(0,0)[cc]{$v$}}
\put(33.721,4.773){\makebox(0,0)[cc]{$e'$}}
\put(31.555,5.613){\line(1,0){5.922}}
\end{picture}


(a) $G$ \hspace{3.8 cm} (b) $G_{v\triangleleft  E_0}$

\caption{Graphs $G$ and $G_{v\triangleleft E_0}$,
where $E_0=\{e_1,e_2,\cdots,e_s\}$}
\relabel{f8}
\end{figure}
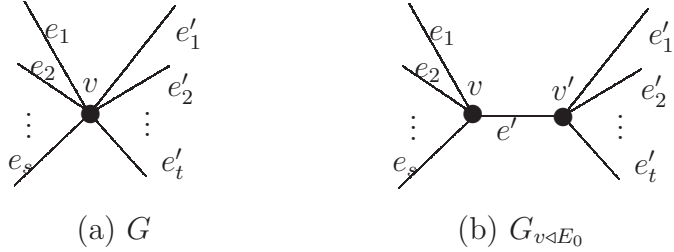

For any subgraph $G_0$ of $G$, 
let $G\diamond G_0$ be the graph below:
$$
G\diamond G_0=(\cdots ((G_{v_1\triangleleft E_1})_{v_2\triangleleft E_2})\cdots )_{v_r\triangleleft E_r},
$$
where $v_1,v_2,\cdots,v_r$ are those vertices in $G_0$ 
with $E_G(v_i)\ne E_{G_0}(v_i)$ and 
$E_i=E_{G_0}(v_i)$.
Clearly, $G_0$ is the subgraph of $G\diamond G_0$ induced 
by $V(G_0)$
and the edges in $E(G\diamond G_0)-E(G)$ 
form a matching of $G\diamond G_0$.
An example is shown in Figure~\ref{f9},
where $G_0=G[\{v_1,v_2,v_3,v_4,v_5,v_6\}]
-\{v_1v_6,v_4v_5\}$
and the new edges in $G\diamond G_0$ are expressed by lines.
If $G\cong K_5$ and $G_0$ is a $5$-cycle, then $G\diamond G_0$ 
is the Petersen graph.

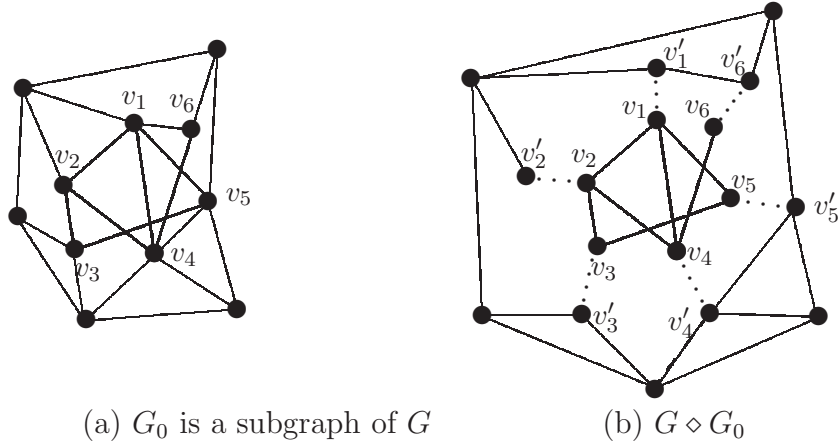
\begin{figure}[h!]
 \centering
\unitlength 2mm 
\linethickness{0.8pt}
\ifx\plotpoint\undefined\newsavebox{\plotpoint}\fi 
\begin{picture}(55.08,26.186)(0,0)
\put(10.282,9.417){\circle*{1.27}}
\put(45.033,9.605){\circle*{1.27}}
\put(5.005,9.715){\circle*{1.27}}
\put(39.756,9.902){\circle*{1.27}}
\put(10.282,9.417){\circle*{1.27}}
\put(45.033,9.605){\circle*{1.27}}
\put(13.85,12.911){\circle*{1.27}}
\put(48.6,13.098){\circle*{1.27}}
\put(12.735,17.668){\circle*{1.27}}
\put(47.486,17.855){\circle*{1.27}}
\put(8.945,18.114){\circle*{1.27}}
\put(43.695,18.301){\circle*{1.27}}
\put(43.695,21.739){\circle*{1.27}}
\put(49.882,20.864){\circle*{1.27}}
\put(52.945,12.489){\circle*{1.27}}
\put(47.195,5.489){\circle*{1.27}}
\put(38.695,5.426){\circle*{1.27}}
\put(35.007,14.551){\circle*{1.27}}
\put(43.57,.426){\circle*{1.27}}
\put(32.07,5.364){\circle*{1.27}}
\put(31.32,21.114){\circle*{1.27}}
\put(51.445,25.551){\circle*{1.27}}
\put(54.445,5.239){\circle*{1.27}}
\put(4.262,13.951){\circle*{1.27}}
\put(39.012,14.139){\circle*{1.27}}
\put(1.574,20.451){\circle*{1.27}}
\put(1.199,11.889){\circle*{1.27}}
\put(5.762,5.076){\circle*{1.27}}
\put(15.762,5.701){\circle*{1.27}}
\put(14.45,23.014){\circle*{1.27}}
\multiput(9.168,18.188)(-.038349206,-.033619048){126}{\line(-1,0){.038349206}}
\multiput(43.918,18.376)(-.03834127,-.033626984){126}{\line(-1,0){.03834127}}
\multiput(4.336,13.952)(.03345,-.223){20}{\line(0,-1){.223}}
\multiput(39.087,14.139)(.03345,-.223){20}{\line(0,-1){.223}}
\multiput(9.093,18.039)(.033648649,-.034648649){148}{\line(0,-1){.034648649}}
\multiput(43.844,18.227)(.033641892,-.034655405){148}{\line(0,-1){.034655405}}
\multiput(4.411,13.877)(.043488889,-.033585185){135}{\line(1,0){.043488889}}
\multiput(39.161,14.065)(.043496296,-.033585185){135}{\line(1,0){.043496296}}
\multiput(10.282,9.343)(-.03323684,.23081579){38}{\line(0,1){.23081579}}
\multiput(45.033,9.531)(-.03326316,.23078947){38}{\line(0,1){.23078947}}
\multiput(4.931,9.641)(.095452632,.033642105){95}{\line(1,0){.095452632}}
\multiput(39.681,9.828)(.095452632,.033642105){95}{\line(1,0){.095452632}}
\multiput(12.884,17.668)(-.03335897,-.10482051){78}{\line(0,-1){.10482051}}
\multiput(47.634,17.855)(-.03334615,-.10482051){78}{\line(0,-1){.10482051}}
\put(9,19.625){\makebox(0,0)[cc]{{\small $v_1$}}}
\put(42.25,19){\makebox(0,0)[cc]{{\small $v_1$}}}
\put(4.563,15.688){\makebox(0,0)[cc]{{\small $v_2$}}}
\put(38.875,15.75){\makebox(0,0)[cc]{{\small $v_2$}}}
\put(5.875,8.188){\makebox(0,0)[cc]{{\small $v_3$}}}
\put(40.188,8.25){\makebox(0,0)[cc]{{\small $v_3$}}}
\put(12.25,9.188){\makebox(0,0)[cc]{{\small $v_4$}}}
\put(46.563,9.25){\makebox(0,0)[cc]{{\small $v_4$}}}
\put(15.938,13.188){\makebox(0,0)[cc]{{\small $v_5$}}}
\put(49.501,14.063){\makebox(0,0)[cc]{{\small $v_5$}}}
\put(46.438,19.063){\makebox(0,0)[cc]{{\small $v_6$}}}
\put(12.171,19.431){\makebox(0,0)[cc]{{\small $v_6$}}}
\thicklines
\multiput(43.368,.305)(.0334821,.0446429){16}{\line(0,1){.0446429}}
\multiput(44.439,1.733)(.0334821,.0446429){16}{\line(0,1){.0446429}}
\multiput(45.511,3.162)(.0334821,.0446429){16}{\line(0,1){.0446429}}
\multiput(46.582,4.59)(.0334821,.0446429){16}{\line(0,1){.0446429}}
\put(45,22.75){\makebox(0,0)[cc]{{\small $v'_1$}}}
\put(48.813,22.313){\makebox(0,0)[cc]{{\small $v'_6$}}}
\put(55.001,12.5){\makebox(0,0)[cc]{{\small $v'_5$}}}
\put(45.313,4.813){\makebox(0,0)[cc]{{\small $v'_4$}}}
\put(40.313,5.188){\makebox(0,0)[cc]{{\small $v'_3$}}}
\thinlines
\multiput(43.626,21.813)(.18382353,-.03308824){34}{\line(1,0){.18382353}}
\multiput(49.876,20.688)(.03333333,.11111111){45}{\line(0,1){.11111111}}
\multiput(51.376,25.688)(.0334186,-.30960465){43}{\line(0,-1){.30960465}}
\multiput(52.813,12.375)(-.033682635,-.041916168){167}{\line(0,-1){.041916168}}
\multiput(54.313,5.125)(-.03316327,.14285714){49}{\line(0,1){.14285714}}
\multiput(51.348,25.595)(-.145101449,-.033514493){138}{\line(-1,0){.145101449}}
\multiput(32.112,5.203)(.081255319,-.033546099){141}{\line(1,0){.081255319}}
\multiput(43.569,.473)(.077892086,.033654676){139}{\line(1,0){.077892086}}
\multiput(47.091,5.413)(-.033542857,-.047552381){105}{\line(0,-1){.047552381}}
\multiput(43.569,.42)(-.033703448,.033710345){145}{\line(0,1){.033710345}}
\put(38.682,5.308){\line(-1,0){6.517}}
\multiput(31.271,21.075)(.6528421,.0332105){19}{\line(1,0){.6528421}}
\multiput(35.055,14.558)(-.033445455,.059245455){110}{\line(0,1){.059245455}}
\multiput(1.629,20.497)(.1732973,.03337838){74}{\line(1,0){.1732973}}
\multiput(14.453,22.967)(-.0328125,-.6306875){16}{\line(0,-1){.6306875}}
\multiput(12.771,17.606)(.03364,.10618){50}{\line(0,1){.10618}}
\multiput(12.666,17.659)(-.283,.0323077){13}{\line(-1,0){.283}}
\multiput(8.987,18.079)(-.10365278,.03358333){72}{\line(-1,0){.10365278}}
\multiput(1.524,20.497)(-.0323077,-.659){13}{\line(0,-1){.659}}
\multiput(1.472,20.497)(.03360465,-.07822093){86}{\line(0,-1){.07822093}}
\multiput(5.834,5.098)(.035604839,.033483871){124}{\line(1,0){.035604839}}
\multiput(10.249,9.25)(.034207547,.033716981){106}{\line(1,0){.034207547}}
\multiput(13.875,12.824)(.03353448,-.12324138){58}{\line(0,-1){.12324138}}
\multiput(15.82,5.676)(-.052066038,.033716981){106}{\line(-1,0){.052066038}}
\multiput(15.662,5.676)(-.4539091,-.0334545){22}{\line(-1,0){.4539091}}
\multiput(4.888,9.67)(.0334545,-.215){22}{\line(0,-1){.215}}
\multiput(5.624,4.94)(-.033731343,.051380597){134}{\line(0,1){.051380597}}
\multiput(1.104,11.825)(.0588209,-.03373134){67}{\line(1,0){.0588209}}
\multiput(32.06,5.413)(-.0325714,.7408095){21}{\line(0,1){.7408095}}
\multiput(47.143,5.413)(1.021143,-.03){7}{\line(1,0){1.021143}}
\multiput(34.88,14.383)(.9985,-.06575){5}{\makebox(0,0)[cc]{$\cdot$}}
\multiput(38.664,5.396)(.1788,.8724){6}{\makebox(0,0)[cc]{$\cdot$}}
\multiput(47.126,5.448)(-.41,.7674){6}{\makebox(0,0)[cc]{$\cdot$}}
\multiput(52.749,12.386)(-.9985,.13125){5}{\makebox(0,0)[cc]{$\cdot$}}
\multiput(47.494,17.747)(.4414,.578){6}{\makebox(0,0)[cc]{$\cdot$}}
\multiput(43.71,18.272)(-.02625,.841){5}{\makebox(0,0)[cc]{$\cdot$}}
\put(35.528,16.082){\makebox(0,0)[cc]{{\small $v'_2$}}}
\end{picture}


(a) $G_0$ is a subgraph of $G$ \hspace{2 cm} 
(b) $G\diamond G_0$

\caption{Graphs $G$ and $G\diamond G_0$}
\relabel{f9}
\end{figure}

Note that $E(G\diamond G_0)-E(G)$ is a matching of
$G\diamond G_0$.
Since $G$ is actually the graph
obtained from  $G\diamond G_0$ by 
contracting all edges in $E(G\diamond G_0)-E(G)$,
applying Lemma~\ref{le2-01}~\ref{le2-1-n4} repeatedly 
on all edges in $E(G\diamond G_0)-E(G)$ implies 
the following result.

\begin{lem}\relabel{le5-3}
Let $G_0$ be a subgraph of $G$ and 
$M=E(G\diamond G_0)-E(G)$.
Then,  for any $N\subseteq E(G)$, 
\begin{equation}\relabel{le5-3-eq1}
\tau_{G}(N)=\tau_{G\diamond G_0}(M\cup N).
\end{equation}
\end{lem}

For a family ${\cal S}=\{E_1,E_2,\cdots,E_k\}$ 
of pairwise disjoint subsets of $E(G)$, 
let $G\diamond {\cal S}$ denote the following graph
obtained by a sequence of $\diamond$-operations
on subgraphs $G[E_1], G[E_2],\cdots, G[E_k]$:
\begin{equation}\relabel{le5-3-eq2}
(\cdots ((G\diamond G[E_1])\diamond G[E_2])\cdots ) \diamond  G[E_k].
\end{equation}
Note that $G\diamond {\cal S}$ is irrelevant to the order of 
$E_1,E_2,\cdots,E_k$
and $(G\diamond {\cal S})/W\cong G$,
where $W=E(G\diamond {\cal S})-E(G)$.
An example of $G\diamond {\cal S}$  
is shown in Figure~\ref{f11}, where 
${\cal S}=\{E_1,E_2,E_3,E_4\}$,
$E_1=E(G[v_1,v_5,v_6])$,  $E_2=E(G[v_2,v_5,v_7,v_9])$, $E_3=E(G[v_3,v_7,v_8])$ and $E_4=E(G[v_4,v_6,v_8,v_9])$.

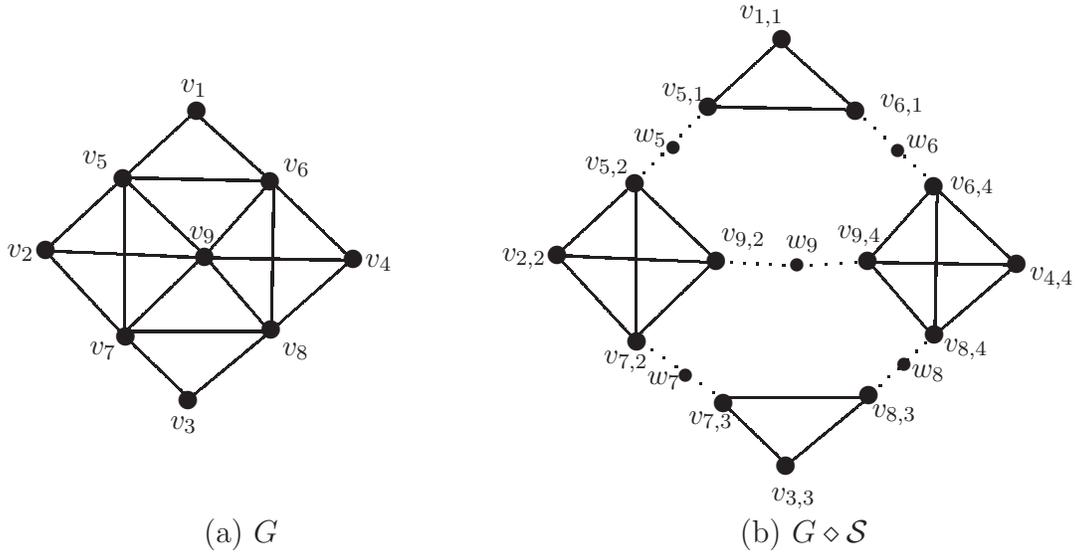
\begin{figure}[h!]
 \centering

\unitlength 2mm 
\linethickness{0.4pt}
\ifx\plotpoint\undefined\newsavebox{\plotpoint}\fi 
\begin{picture}(68.599,32.433)(0,0)
\put(22.149,16.327){\circle*{1.27}}
\put(66.244,16.012){\circle*{1.27}}
\put(11.743,26.208){\circle*{1.27}}
\put(50.635,30.938){\circle*{1.27}}
\put(1.705,16.905){\circle*{1.27}}
\put(35.709,16.59){\circle*{1.27}}
\put(11.165,6.972){\circle*{1.27}}
\put(50.898,2.557){\circle*{1.27}}
\put(16.683,11.649){\circle*{1.27}}
\put(56.416,7.235){\circle*{1.27}}
\put(60.778,11.334){\circle*{1.27}}
\put(12.269,16.485){\circle*{1.27}}
\put(46.273,16.169){\circle*{1.27}}
\put(56.364,16.169){\circle*{1.27}}
\put(16.631,21.478){\circle*{1.27}}
\put(55.523,26.208){\circle*{1.27}}
\put(60.726,21.162){\circle*{1.27}}
\put(6.855,21.688){\circle*{1.27}}
\put(45.747,26.418){\circle*{1.27}}
\put(40.859,21.372){\circle*{1.27}}
\put(7.013,11.176){\circle*{1.27}}
\put(46.746,6.762){\circle*{1.27}}
\put(41.017,10.861){\circle*{1.27}}
\thicklines
\multiput(6.938,21.706)(.7519231,-.0161538){13}{\line(1,0){.7519231}}
\multiput(45.83,26.436)(.7519231,-.0161538){13}{\line(1,0){.7519231}}
\multiput(16.713,21.496)(-.0168473282,-.0188549618){262}{\line(0,-1){.0188549618}}
\multiput(60.808,21.181)(-.0168473282,-.0188587786){262}{\line(0,-1){.0188587786}}
\multiput(12.299,16.556)(-.0176930693,.0168250825){303}{\line(-1,0){.0176930693}}
\multiput(46.303,16.24)(-.0176930693,.0168250825){303}{\line(-1,0){.0176930693}}
\put(6.938,21.654){\line(0,-1){10.144}}
\put(40.942,21.338){\line(0,-1){10.143}}
\multiput(7.043,11.51)(.0173433333,.01682){300}{\line(1,0){.0173433333}}
\multiput(41.047,11.195)(.0173433333,.0168166667){300}{\line(1,0){.0173433333}}
\multiput(12.246,16.556)(.0168416988,-.0194826255){259}{\line(0,-1){.0194826255}}
\multiput(56.341,16.24)(.0168416988,-.0194787645){259}{\line(0,-1){.0194787645}}
\put(16.608,11.51){\line(-1,0){9.565}}
\put(56.341,7.096){\line(-1,0){9.565}}
\multiput(16.871,21.391)(-.0158,-.9828){10}{\line(0,-1){.9828}}
\multiput(60.966,21.076)(-.0158,-.9828){10}{\line(0,-1){.9828}}
\put(11.563,27.698){\makebox(0,0)[cc]{{\small $v_1$}}}
\put(.053,16.766){\makebox(0,0)[cc]{{\small $v_2$}}}
\put(10.88,5.309){\makebox(0,0)[cc]{{\small $v_3$}}}
\put(23.861,16.083){\makebox(0,0)[cc]{{\small $v_4$}}}
\put(5.046,22.915){\makebox(0,0)[cc]{{\small $v_5$}}}
\put(18.395,22.39){\makebox(0,0)[cc]{{\small $v_6$}}}
\put(5.519,10.249){\makebox(0,0)[cc]{{\small $v_7$}}}
\put(18.343,10.196){\makebox(0,0)[cc]{{\small $v_8$}}}
\put(12.141,17.922){\makebox(0,0)[cc]{{\small $v_9$}}}
\multiput(22.284,16.24)(-.9985,.0158){10}{\line(-1,0){.9985}}
\multiput(66.379,15.925)(-.9985,.0158){10}{\line(-1,0){.9985}}
\multiput(12.299,16.398)(-.3678966,.0163103){29}{\line(-1,0){.3678966}}
\multiput(46.303,16.083)(-.3678966,.0163103){29}{\line(-1,0){.3678966}}
\put(1.63,16.871){\line(0,1){0}}
\multiput(60.651,21.233)(.0178507937,-.0168507937){315}{\line(1,0){.0178507937}}
\multiput(66.274,15.925)(-.0196370107,-.0168327402){281}{\line(-1,0){.0196370107}}
\multiput(11.142,6.938)(.0203955224,.0168656716){268}{\line(1,0){.0203955224}}
\multiput(50.875,2.523)(.0203955224,.0168656716){268}{\line(1,0){.0203955224}}
\multiput(16.608,11.458)(.0190452962,.0168466899){287}{\line(1,0){.0190452962}}
\put(22.074,16.293){\line(0,1){0}}
\multiput(22.074,16.293)(-.0173464052,.0168333333){306}{\line(-1,0){.0173464052}}
\multiput(16.766,21.444)(-.0187309091,.0168145455){275}{\line(-1,0){.0187309091}}
\multiput(55.658,26.174)(-.0187309091,.0168181818){275}{\line(-1,0){.0187309091}}
\multiput(11.615,26.068)(-.0182519084,-.0168473282){262}{\line(-1,0){.0182519084}}
\multiput(50.507,30.799)(-.0182519084,-.016851145){262}{\line(-1,0){.0182519084}}
\multiput(6.833,21.654)(-.0181373239,-.0168415493){284}{\line(-1,0){.0181373239}}
\multiput(40.837,21.338)(-.0181373239,-.0168380282){284}{\line(-1,0){.0181373239}}
\multiput(1.682,16.871)(.0168461538,-.0181923077){312}{\line(0,-1){.0181923077}}
\multiput(35.686,16.556)(.0168461538,-.0181923077){312}{\line(0,-1){.0181923077}}
\multiput(6.938,11.195)(.017448,-.01682){250}{\line(1,0){.017448}}
\multiput(46.671,6.78)(.017448,-.016816){250}{\line(1,0){.017448}}
\multiput(40.802,21.356)(.611,.6175){9}{{\rule{.8pt}{.8pt}}}
\multiput(55.465,26.191)(.657,-.62412){9}{{\rule{.8pt}{.8pt}}}
\multiput(41.117,10.792)(.81843,-.57057){8}{{\rule{.8pt}{.8pt}}}
\multiput(56.516,7.271)(.60071,.56314){8}{{\rule{.8pt}{.8pt}}}
\put(43.412,23.756){\circle*{.841}}
\put(58.338,23.546){\circle*{.841}}
\put(51.663,15.925){\circle*{.841}}
\put(44.253,8.62){\circle*{.841}}
\put(58.759,9.303){\circle*{.841}}
\put(49.301,32.433){\makebox(0,0)[cc]{{\small $v_{1,1}$}}}
\put(44.149,27.283){\makebox(0,0)[cc]{{\small $v_{5,1}$}}}
\put(58.642,26.54){\makebox(0,0)[cc]{{\small $v_{6,1}$}}}
\put(39.02,22.452){\makebox(0,0)[cc]{{\small $v_{5,2}$}}}
\put(33.52,16.357){\makebox(0,0)[cc]{{\small $v_{2,2}$}}}
\put(40.284,9.519){\makebox(0,0)[cc]{{\small $v_{7,2}$}}}
\put(48.088,17.621){\makebox(0,0)[cc]{{\small $v_{9,2}$}}}
\put(45.933,5.506){\makebox(0,0)[cc]{{\small $v_{7,3}$}}}
\put(51.433,.377){\makebox(0,0)[cc]{{\small $v_{3,3}$}}}
\put(58.122,5.952){\makebox(0,0)[cc]{{\small $v_{8,3}$}}}
\put(62.879,10.337){\makebox(0,0)[cc]{{\small $v_{8,4}$}}}
\put(68.599,15.058){\makebox(0,0)[cc]{{\small $v_{4,4}$}}}
\put(63.301,21.095){\makebox(0,0)[cc]{{\small $v_{6,4}$}}}
\put(55.818,17.844){\makebox(0,0)[cc]{{\small $v_{9,4}$}}}
\put(41.909,24.386){\makebox(0,0)[cc]{{\small $w_5$}}}
\put(60.034,23.996){\makebox(0,0)[cc]{{\small $w_6$}}}
\put(52.027,17.324){\makebox(0,0)[cc]{{\small $w_9$}}}
\put(42.821,8.228){\makebox(0,0)[cc]{{\small $w_7$}}}
\put(60.392,8.729){\makebox(0,0)[cc]{{\small $w_8$}}}
\multiput(56.463,16.047)(-.9776,-.0524){6}{{\rule{.8pt}{.8pt}}}
\multiput(51.575,15.785)(-.86717,.04367){7}{{\rule{.8pt}{.8pt}}}
\end{picture}


(a)  $G$\hspace{6 cm}
(b) $G\diamond {\cal S}$

\caption{Graphs $G$ and $G\diamond {\cal S}$}
\relabel{f11}
\end{figure}

Assume that $V(G)=\{v_1,v_2,\cdots,v_n\}$.
By the above definition,  
if  ${\cal S}=\{E_1,E_2,\cdots,E_k\}$ 
is a partition of $E(G)$, then 
$G\diamond {\cal S}$ is actually the graph with vertex set:
\begin{equation}\relabel{le5-3-eq3}
\{w_i: v_i\in V'\}\cup
\bigcup_{j=1}^k V_j,
\end{equation}
where $V'=V(G)-\{v_i\in V(G): \exists j, N_G(v_i)=E_j\}$
and 
$V_j=\{ v_{i,j}: v_i\in V(G[E_j])\}$, and  edge set: 
\begin{equation}\relabel{le5-3-eq4}
\bigcup_{j=1}^k E'_j\cup 
\bigcup_{v_i\in V'}
\left \{w_iv_{i,j}: E_G(v_i)\cap E_j\ne \emptyset,1\le j\le k\right \},
\end{equation}
where each $E'_j$ is a copy of $E_j$ by changing 
the ends $v_{s}$ and $v_{t}$ of each edge $e$ in $E_j$ 
to $v_{s,j}$ and $v_{t,j}$.
As each edge in $G$ has exactly one copy in 
$G\diamond {\cal S}$,
$E(G\diamond {\cal S})$ is also considered 
as the union of 
$E(G)$ and $\bigcup_{v_i\in V'}
\{w_iv_{i,j}: E_G(v_i)\cap E_j\ne \emptyset,1\le j\le k\}$.
An example for the labels of vertices and edges 
in $G\diamond {\cal S}$ is given in Figure~\ref{f11}.

Some basic facts on $G\diamond {\cal S}$ 
follow directly. 

\begin{lem}\relabel{le5-30}
Let $M=E(G\diamond {\cal S})-E(G)$.
Then 
\begin{enumerate}
\item $|V(G\diamond {\cal S})-V(G)|=|V'|$;
\item the component number 
of $(G\diamond {\cal S})[M]$ is equal to $|V'|$;
\item each component of $(G\diamond {\cal S})[M]$
is a star $S_j$ with its center $w_j\in V(G\diamond {\cal S})-V(G)$
and its size equal to the number of different 
sets $E_i$ with $E_G(v_j)\cap E_i\ne \emptyset$;
\item $\{w_i: v_i\in V'\}$ is an independent set 
in $G\diamond {\cal S}$ and its removal from $G\diamond {\cal S}$
results in $k$ components isomorphic to
$G[E_1], G[E_2],\cdots, G[E_k]$ respectively;
\item  $G\diamond {\cal S}[V_j] \cong G[E_j]$
for each $j$;
\item 
$E_{G\diamond {\cal S}}(V_{j_1},V_{j_2})=\emptyset$
for all $1\le j_1<j_2\le k$.
\end{enumerate}
\end{lem}

Applying Lemma~\ref{le5-3}, we get the following conclusion.

\begin{lem}\relabel{le5-4}
For any partition ${\cal S}=\{E_1,E_2,\cdots,E_k\}$ of
$E(G)$ and any $N\subseteq E(G)$, we have 
\begin{equation}\relabel{le5-4-eq1}
\tau_G(N)=\tau_{G\diamond {\cal S}}(M\cup N),
\end{equation}
where $M=E(G\diamond {\cal S})-E(G)$.
\end{lem}

If $G[E_i]$ is a complete graph in $G$ 
for all $i=1,2,\cdots,k$, 
applying Lemma~\ref{le5-4} and Theorem~\ref{th4-2}
gets the following expression on $\tau_G$. 

\begin{theo}\relabel{th5-4}
Assume that ${\cal S}=\{E_1,E_2,\cdots,E_k\}$ 
is a partition of $E(G)$ such that $G[E_i]$ is 
a complete graph for all $i=1,2,\cdots,k$.
Then  
\begin{equation}\relabel{th5-4-eq1}
\tau_G=\prod_{i=1}^k  n_i^{n_i-2-|n'_i|}
\sum_{T\in \setst_{(G\diamond {\cal S})/E(G)}}
\prod_{i=1}^k 
n_i^{|E_T(v'_i)|},
\end{equation}
where $n_i=|V(G[E_i])|$, 
$n'_i$ is the size of the set 
$\{v_s\in V(G[E_i]):  E_G(v_s)\not\subseteq E_i\}$ 
and $v'_i$ is the new vertex in $(G\diamond {\cal S})/E(G)$
produced by contracting $V_i$ in $G\diamond \sets$.
\end{theo}

Note that $(G\diamond {\cal S})/E(G)$ is the bipartite graph 
with a bipartition $\{w_i: v_i\in V'\}$ 
and $\{v'_j:1\le j\le k\}$ and edge set 
$\{w_iv'_j: E_G(v_i)\cap E_j\ne \emptyset\}$.

\subsection{Application \relabel{sec5-02}}

Obviously Theorem~\ref{th5-4} can be applied to 
the graph in Figure~\ref{f11} (a).
Actually this graph is 
the middle graph of $K_4-e$ (i.e., the graph with one 
edge removed from $K_4$).
For a graph $H$ with vertex set $\{u_1,u_2,\cdots,u_k\}$,
the {\it middle graph}\relabel{middleg} of $H$, 
denoted by $M(H)$,  
is the graph obtained from its line graph $L(H)$ 
and the empty graph $H-E(H)$ 
by adding edges joining each vertex $u_i$ in $H-E(H)$ 
to all those vertices in $L(H)$ 
which correspond to edges in the set $E_H(u_i)$.

Applying generalized Wye-Delta transform and 
Delta-Wye transform,
Yan~\cite{yan0} gave a relation between 
$\setst_{M(H)}$ and $\setst_{S(H)}$,
where $S(H)$ is the graph obtained from $H$ by subdividing 
each edge in $H$ exactly once.
Such a relation actually follows from 
Theorem~\ref{th5-4} directly. 
Observe that the edge set of $M(H)$ 
has a partition 
${\cal S}=\{E_1,E_2,\cdots, E_k\}$, 
where  
each $E_i$ is the set of edges in the subgraph 
of $M(H)$ induced by $\{u_i\}\cup N_{M(H)}(u_i)$.
Clearly each $M(H)[E_i]$ is a complete graph of order 
$d_H(u_i)+1$ and contains exactly $d_H(u_i)$ 
vertices $u$ such $E_{M(H)}(u)\not\subseteq M(H)[E_i]$.
Also note that $M(H)\diamond {\cal S}/E(M(H))$ 
is actually the graph $S(H)$.
Thus, applying Theorem~\ref{th5-4},
we have 
\begin{equation}\relabel{th5-4-eq2}
\tau_{M(H)}
=
\sum_{T\in \setst_{S(H)}}
\prod_{i=1}^k (d_H(u_i)+1)^{|E_T(u_i)|-1}.
\end{equation}
Similarly, a relation between 
$\setst_{L(H)}$ and $\setst_{S(H)}$ can be obtained:
\begin{equation}\relabel{th5-4-eq3}
\tau_{L(H)}
=
\sum_{T\in \setst_{S(H)}}
\prod_{i=1}^k d_H(u_i)^{|E_T(u_i)|}.
\end{equation}
It is not difficult to verify that 
(\ref{eq1-3}) can be obtained from 
(\ref{th5-4-eq3}). 

\resection{
A factorization of $\tau_G(W)$ 
\relabel{sec6}
}

If a simple graph $G=(V,E)$ contains a clique $U$ 
and a partition $S_1, S_2$ of $V-U$ 
with $S_1\cap N_G(S_2)=\emptyset$, 
then the chromatic polynomial $\chi(G,\lambda)$  has 
the following factorization due to Zykov~\cite{zyk}
(see~\cite{dong0,rea1} also):
\begin{equation}\relabel{eq6-001}
\chi(G,\lambda)=\frac{\chi(G[U\cup S_1],\lambda)
\chi(G[U\cup S_2],\lambda)}
{\chi(G[U],\lambda)}.
\end{equation} 
In this section, 
we find a similar expression for $\tau_G(W)$ 
for any $W\subseteq E(G)-E(G[U])$
by applying results in Sections~\ref{sec2}, ~\ref{sec3} and~\ref{sec4}. 

\begin{theo}\relabel{th5-10}
Let $G=(V,E)$ be a connected and loopless multigraph
and $U$ be a clique of $G$.
If $S_1$ and $S_2$ form a partition 
of $V-U$ with 
$N_G[S_1]\cap N_G[S_2]=\emptyset$, 
as shown in Figure~\ref{f7}, 
then, for any $W\subseteq E-E(G[U])$,
\begin{equation}\relabel{eq5-10}
\tau_{G}(W) 
=\frac{\tau_{G[U\cup S_1]}(W_1)
\cdot 
\tau_{G[U\cup S_2]}(W_2)}
{|U|^{|U|-2}},
\end{equation}
where $W_i=W\cap E(G[U\cup S_i])$.
\end{theo}

\begin{figure}[htbp]
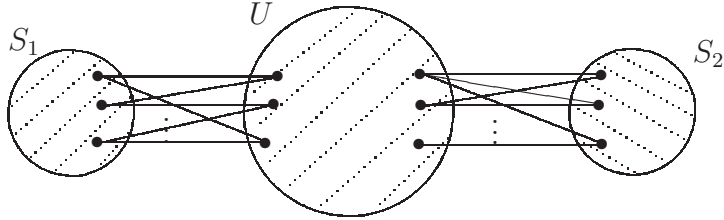

 \centering
\unitlength .99mm 
\linethickness{0.4pt}
\ifx\plotpoint\undefined\newsavebox{\plotpoint}\fi 



\caption{A clique $U$ and a partition $S_1,S_2$ of $V-U$ with $N_G[S_1]\cap N_G[S_2]=\emptyset$}
\relabel{f7}
\end{figure}

\proof 
Let $M=E_G(U)$, $R=W\cap M$, $G_i=G[U\cup S_i]$,
$M_i=M\cap E(G_i)$, 
$R_i=R\cap E(G_i)$ and $N_i=(W-R)\cap E(G_i)$ for $i=1,2$.
Thus $W_i$ is the disjoint union of $R_i$ and $N_i$.
We will prove this result by the following claims.

\inclaim 
(\ref{eq5-10}) holds if each component of $G[M]$ 
is a star with a center in $S_1\cup S_2$.

Let $u_i$ represent the vertex in $G_i\bullet U$ 
after identifying all vertices in $U$ as one vertex.
For $T_i\in \setst_{G_i\bullet U}(N_i)$ 
for $i=1,2$, 
let $T_1\cdot T_2$ denote the tree
obtained from $T_1$ and $T_2$ by identifying $u_1$ 
and $u_2$ as one vertex. 
By the definition of $G\bullet U$ and the given condition, 
we have 
$$
\setst_{G\bullet U}(N_1\cup N_2)
=\{T_1\cdot T_2: T_i\in \setst_{G_i\bullet U}(N_i),i=1,2\}.
$$
Thus, for $i=1,2$, by Theorem~\ref{th5-3} with $k=1$, 
\begin{eqnarray}\relabel{eq5-11-1}
\tau_{G_i}(W_i)
&=&|U|^{|U|-2}
\sum_{T_i\in \setst_{G_i\bullet U}(N_i)}
|U|^{-|M_i-E(T_i)|} 
(1+|U|)^{|(M_i-R_i)-E(T_i)|}
\end{eqnarray}
 and  
\begin{eqnarray}\relabel{eq5-11}
\tau_G(W) &=&
|U|^{|U|-2}
\sum_{T\in \setst_{G\bullet U}(N_1\cup N_2)}
|U|^{-|M- E(T)|}
(1+|U|)^{|(M-(R_1\cup R_2))-E(T)|}
\nonumber \\
&=&|U|^{|U|-2}
\prod_{i=1}^2  
\sum_{T_i\in \setst_{G_i\bullet U}(N_i)}
|U|^{-|M_i-E(T_i)|}
(1+|U|)^{|(M_i-R_i)-E(T_i)|}
\nonumber \\
&=&\frac{\tau_{G_1}(W_1) \tau_{G_2}(W_2)}
{|U|^{|U|-2}}.
\end{eqnarray}
Thus Claim~\thecountclaim\   holds.

\inclaim (\ref{eq5-10}) holds if $R=M$ and $G[M]$ is a forest.

Let $G'=G\star W-M$ and $W'=E(G\star W)-E(G)$.
For $i=1,2$, 
let $W_i=W\cap E(G_i)$,
$G_i'=G_i\star W_i-R_i$
and 
$W_i'=E(G\star W_i)-E(G_i)$.
By 
Lemma~\ref{le2-3}~\ref{le2-3-n1},
\begin{equation}\relabel{eq6-01}
\tau_{G}(W)
=\tau_{G'}(W')|.
\end{equation}
and 
\begin{equation}\relabel{eq6-02}
\tau_{G_i}(W_i)
=|\setst_{G'_i}(W'_i)|,
\quad i=1,2.
\end{equation}
Note that $U$ is a clique of $G'$ 
and $G'[E_{G'}(U)]$ is a star.
By Claim 1, 
\begin{equation}\relabel{eq6-03}
\tau_{G'}(W')|=
\frac{\tau_{G'_1}(W'_1) \tau_{G'_2}(W'_2)}
{|U|^{|U|-2}}.
\end{equation}
Thus, Claim 2 follows from 
(\ref{eq6-01}), (\ref{eq6-02}) and (\ref{eq6-03}).

\inclaim (\ref{eq5-10}) holds.

For any $M'\subseteq M$, 
$G[M']$ is not a forest \iff  
$G_i[M'_i]$ is not a forest for some $i=1,2$,
where $M'_i=M'\cap E_G(U)$.
Thus 
\begin{equation}\relabel{eq6-04}
\tau_G(M'\cap N)=
\frac{\tau_G(M'_1\cap N_1)
\tau_G(M'_2\cap N_2)}{|U|^{|U|-2}}=0.
\end{equation}
Let $R$ be a fixed subset of $M$ such that 
$G[R]$ is a forest.
For $i=1,2$,
let $\R_i=\{R'_i:R_i\subseteq R_i'\subseteq M\cap E(G_i), 
G[R_i'] \mbox{ is a forest}\}$,
where $R_i=R\cap E(G_i)$.
Note that 
\begin{equation}\relabel{eq6-05}
\setst_{G}(R\cup N)=
\bigcup_{R_i'\in \R_i
\atop i=1,2}
\setst_{G-(M-(R_1'\cup R'_2))}
(R_1'\cup R_2'\cup N_1\cup N_2),
\end{equation}
where for distinct order pairs 
$(R'_1,R'_2)$ and $(R''_1,R''_2)$ in the above union,
the two corresponding sets 
$\setst_{G-(M-(R_1'\cup R'_2))}
(R_1'\cup R_2'\cup N_1\cup N_2)$ 
and $\setst_{G-(M-(R_1''\cup R''_2))}
(R_1''\cup R_2''\cup N_1\cup N_2)$ 
are disjoint. 
Similarly,
\begin{equation}\relabel{eq6-06}
\setst_{G_i}(R_i\cup N_i)=
\bigcup_{R'_i\in \R_i}
\setst_{G_i-(M\cap E(G_i)-R'_i)}(R'_i\cup N_i),
\quad i=1,2,
\end{equation}
where the above union is disjoint union for both $i=1,2$.
By Claim 2, for any $R'_i\in \R_i$ for $i=1,2$, we
have 
\begin{equation}\relabel{eq6-07}
\tau_{G-(M-(R_1'\cup R'_2))}(R_1'\cup R_2'\cup N_1\cup N_2)
=\frac{\prod_{1\le i\le 2} 
\tau_{G_i-(M\cap E(G_i)-R'_i)}(R'_i\cup N_i)}{|U|^{|U|-2}}.
\end{equation}

Thus, Claim 3 follows from 
(\ref{eq6-05}), (\ref{eq6-06}) and (\ref{eq6-07}),
and the result is proven.
\proofend

\begin{cor}\relabel{co5-1}
Let $G=(V,E)$ be any connected multigraph 
and $U$ be a clique of $G$.
If $w$ is a vertex in $V-U$ with 
$N_G[w]\cap N_G[V-(U\cup \{w\})]=\emptyset$,
then 
\begin{equation}\relabel{eq5-8}
\tau_{G}
=\tau_{G-w}\left ( d(w)(1+1/|U|)^{d(w)-1}\right ).
\end{equation}
\end{cor}

\proof Let $S_1=\{w\}$ and $S_2=V-\{w\}$.
As $N_G[w]\cap N_G[V-(U\cup \{w\})]=\emptyset$,
by applying Theorem~\ref{th5-10},
\begin{equation}\relabel{eq6-08}
\tau_G=\frac{\tau_{G_1}\cdot \tau_{G_2}}
{|U|^{|U|-2}},
\end{equation}
where $G_1=G[U\cup \{w\}]$ and $G_2=G[U\cup S_2]=G-w$.
By Theorem~\ref{th5-3},  
\begin{equation}\relabel{eq6-09}
\tau_{G_1}=d(w)\cdot |U|^{|U|-2}\cdot (1+1/|U|)^{d(w)-1}.
\end{equation}
The result then follows from (\ref{eq6-08})
and (\ref{eq6-09}). 
\proofend

\noindent {\bf Remarks}:
(a) The condition for (\ref{eq6-001})
is weaker than the one for (\ref{eq5-10}),
as (\ref{eq5-10}) holds with an extra condition 
$N_G[S_1]\cap N_G[N_2]=\emptyset$.

(b)
Note that when $W=\emptyset$, (\ref{eq5-10}) 
is equivalent to the following equality 
\begin{equation}\relabel{eq6-10}
T_G(x,y)T_{K_{|U|}}(x,y)=T_{G[U\cup S_1]}(x,y)
T_{G[U\cup S_2]}(x,y)
\end{equation}
when $(x,y)=(1,1)$, 
where $T_G(x,y)$ is the Tutte polynomial of $G$.

(c) Under the condition 
that $S_1\cap N_G(S_2)=\emptyset$, 
(\ref{eq6-001}) implies that 
(\ref{eq6-10}) holds when $y=0$,
as $T_G(1-x,0)=x^{-c(G)}(-1)^{|V|-c(G)}\chi(G,x)$ holds 
for any simple graph $G$ (see \cite{bry, ellis}),
where $c(G)$ is the number of components of $G$.
Furthermore,  
(\ref{eq6-10}) also holds 
for graph $G$ satisfying condition 
$N_G[S_1]\cap N_G[N_2]=\emptyset$
when $(x,y)=(2,2)$, 
as $T_G(2,2)=2^{|E(G)|}$ holds for any graph $G$.
We have verified (\ref{eq6-10}) for 
some graphs $G$  satisfying the same condition
when $(x,y)=(0,-1)$,
but we are not sure if it holds 
for all graphs $G$
satisfying this condition.

\begin{prob}\relabel{prob6-1}
Let $U$ be a clique of $G=(V,E)$.
If $V-U$ has a partition $S_1$ and $S_2$ 
with $N_G[S_1]\cap N_G[N_2]=\emptyset$,  
does (\ref{eq6-10}) hold at $(x,y)=(0,-1)$?
\end{prob}

\end{document}